\documentclass[hidelinks, 11pt]{article}

\usepackage[UKenglish]{babel}
\usepackage[utf8]{inputenc}
\usepackage[T1]{fontenc}
\usepackage{erewhon}
\usepackage{amsmath}
\usepackage{mathrsfs}
\usepackage{amsthm}
\usepackage{amsfonts}
\usepackage{amssymb}
\usepackage{mathtools}
\usepackage{dsfont} 
\usepackage{bm} 

\usepackage{quiver}
\usetikzlibrary{nfold}
\usepackage{adjustbox}
\usepackage{bussproofs}

\usepackage[dvipsnames]{xcolor}
\usepackage{hyperref}
\usepackage{fullpage}
\usepackage{microtype}
\usepackage{booktabs}
\usepackage{enumitem} 
\usepackage{comment}
\usepackage{authblk} 

\usepackage[autostyle,italian=guillemets]{csquotes}
\usepackage[backend=biber]{biblatex}
\tikzset{node distance=2cm, auto}
\tikzcdset{row sep/normal=2.7em,column sep/normal=4em}

\allowdisplaybreaks

\setlist[description]{font=\normalfont\bfseries\:\!}

\renewcommand\labelenumi{\normalfont\textbf{[\arabic{enumi}]}}
\renewcommand\theenumi\labelenumi

\hypersetup{
    colorlinks=true,
    linkcolor=gray,
    filecolor=gray,      
    urlcolor=gray,
    citecolor=gray
    }
\numberwithin{equation}{section}

\newtheorem{theorem}{Theorem}[section]
\newtheorem{proposition}[theorem]{Proposition}
\newtheorem{lemma}[theorem]{Lemma}
\newtheorem{corollary}[theorem]{Corollary}

\theoremstyle{definition}
\newtheorem{definition}[theorem]{Definition}
\newtheorem{example}[theorem]{Example}

\newtheorem{remark}[theorem]{Remark}

\newenvironment{contentX}[1][blank]{\noindent{\large\bfseries\scshape\MakeLowercase{#1}}.\thinspace\noindent}{}
\newenvironment{content}[1][blank]{\vskip .5cm\begin{contentX}[#1]}{\end{contentX}}
\newenvironment{acknowledgements}{\begin{content}[Acknowledgements]}{\end{content}}

\newenvironment{conventions}{\begin{content}[Conventions]}{\end{content}}

\renewenvironment{abstract}{\begin{contentX}[Abstract]}{\end{contentX}}

\DeclareSymbolFont{sfoperators}{OT1}{cmss}{m}{n}

\makeatletter
\def\operator@font{\mathgroup\symsfoperators}
\makeatother

\newcommand{\CategoryFont}[1]{\mathsf{\uppercase{#1}}}
\newcommand{\FunctorFont}[1]{\mathsf{#1}}

\renewcommand{\,}{,\dots,}

\renewcommand{\=}{\mathrel{\mathop:}=}
\renewcommand{\epsilon}{\varepsilon}
\renewcommand{\o}{\circ}

\newcommand{\<}{\langle}
\renewcommand{\>}{\rangle}
\newcommand{\op}{\FunctorFont{op}}

\renewcommand{\*}{\ast}
\newcommand{\nto}{\mathrel{\mkern3mu\vcenter{\hbox{$\scriptscriptstyle|\:$}}\mkern-12mu{\to}}}

\newcommand{\ntimes}{\times\dots\times}
\newcommand{\ex}{\mathsf{ex}}

\newcommand{\Real}{\mathbb{\R}}

\newcommand{\Mod}{\CategoryFont{Mod}}

\newcommand{\KL}{\FunctorFont{KL}}
\newcommand{\id}{\FunctorFont{id}}

\newcommand{\X}{\mathbb{X}}

\newcommand{\Set}{\CategoryFont{Set}}

\DeclareMathOperator{\DO}{DO}

\newcommand{\q}{\mathbf{q}}

\renewcommand{\d}{\FunctorFont{d}}
\newcommand{\T}{\mathrm{T}}
\newcommand{\TT}{\mathbb{T}}

\newcommand{\SMan}{\CategoryFont{SMan}}

\newcommand{\Aff}{\CategoryFont{Aff}}

\DeclareMathOperator{\GL}{GL}

\DeclareMathOperator{\Ev}{Ev}
\DeclareMathOperator{\ev}{ev}

\DeclareMathOperator{\Sym}{Sym}

\newcommand{\gl}{\mathfrak{gl}}
\newcommand{\lnr}{\mathscr{L}}

\newcommand{\Simple}{\FunctorFont{S}}

\newcommand{\g}{\mathfrak{g}}

\newcommand{\LieAlg}{\CategoryFont{LieAlg}}

\newcommand{\e}{\FunctorFont{e}}
\newcommand{\m}{\FunctorFont{m}}
\renewcommand{\i}{\FunctorFont{i}}
\newcommand{\Gr}{\CategoryFont{Grp}}
\newcommand{\SQ}[1]{[\![#1]\!]}
\newcommand{\PLIVF}{PLIVF }
\newcommand{\PLIVFs}{PLIVFs }

\newcommand{\A}{{\bm{A}}}
\newcommand{\B}{{\bm{B}}}
\newcommand{\C}{{\bm{C}}}
\newcommand{\R}{\bm{R}}

\newcommand{\multi}{\mathsf{multi}}
\DeclareMathOperator{\Shuffle}{Shuffle}

\title{General linear groups in tangent categories}
\author{\uppercase{Marcello Lanfranchi}}
\affil{\normalsize\textit{Macquarie University, School of Mathematical and Physical Sciences}}
\date{}

\begin{document}

\maketitle

\begin{abstract}
General linear groups are well-studied Lie groups that encode the linear symmetries of a finitely-dimensional linear space. In differential geometry, they are an important ingredient of various constructions, including fibre bundles and the frame bundle of a principal bundle. They also appear in algebraic geometry as affine group schemes. Their Lie algebra, known as the general linear Lie algebra, is isomorphic to the Lie algebra of square matrices, whose Lie bracket is given by the commutator. Building on previous work on Lie groups, in this paper we introduce general linear group objects in the context of tangent categories. We define the general linear group of a differential object as a classifier of the Kleisli parametrized linear automorphisms of the differential object, a class of maps that generalize linear automorphisms in a suitable sense. We prove that such a group can be entirely constructed via a tangent triple equalizer and use this characterization to prove that the Lie algebra of general linear group objects is isomorphic to the Lie algebra of linear endomorphisms.
\end{abstract}

\begin{acknowledgements}
Mathematics is an emergent phenomenon that happens through a complex process of interactions with other researchers and not the achievement of isolated minds. Although the main ideas of this paper have been developed and carried over by its author, they could not exist without the numerous useful discussions and suggestions of various people of the growing community of tangent and differential categories. In particular, the author would like to thank JS Lemay for his support, mentorship and suggestions on this work. For this work, the author was supported by the AFOSR under award number FA9550-24-1-0008.
\end{acknowledgements}

\par\noindent\rule{\textwidth}{1pt}


\tableofcontents

\section{Introduction}
\label{section:introduction}
The $n$-th general linear group is the group $\GL_n(\Real)$ of linear automorphisms of $\Real^n$. Thanks to the differentiability of linear maps, $\GL_n(\Real)$ comes equipped with a smooth structure that makes it into a Lie group. A Lie group is a group object in the category $\SMan$ of finitely-dimensional smooth manifolds. More concretely, a Lie group consists of a group whose underlying set is a smooth manifold and whose multiplication and inverse maps are smooth with respect to that smooth structure. Every Lie group $G$ comes with an associated Lie algebra $\g$, which is classically defined as the Lie algebra of vector fields on $G$ that are invariant under the left action $L_g(h)\=gh$ of the group. The Lie algebra $\g$ is also isomorphic to the tangent space $\T_\e G$ of $G$ at the unit.

In previous work~\cite{lanfranchi:lie-correspondence-tangent-cats}, we internalized this correspondence between Lie groups and Lie algebras in tangent categories.

\begin{content}[Tangent categories]
Introduced by Rosick\'y in~\cite{rosicky:tangent-cats} and later generalized by Cockett and Cruttwell in~\cite{cockett:tangent-cats}, tangent categories offer a minimal categorical context for differential geometry. A tangent category consists of a collection of objects, interpreted as generalized geometric spaces, a collection of morphisms, interpreted as smooth functions, and an assignment, $\T$, which sends each object $M$ to a new object $\T M$ and each morphism $f\colon M\to N$ to a new morphism $\T f\colon\T M\to\T N$.

The object $\T M$ should be regarded as the space of tangent vectors of $M$, known as the \textbf{tangent bundle} of $M$, while $\T f$ is interpreted as the smooth function which sends each tangent vector $u\in\T M$ of $M$ at $x$ to the tangent vector $\d_xf(u)$, where $\d_xf$ denotes the differential of $f$ at $x$. The assignment $\T$ is functorial, in that $\T$ preserves identities and composition of morphisms.

A tangent category also comes with a family of natural transformations
\begin{align*}
&p_M\colon\T M\to M         &&z_M\colon M\to\T M            &&s_M\colon\T_2M\to\T M\\
&l_M\colon\T M\to\T\T M     &&                              &&c_M\colon\T\T M\to\T\T M
\end{align*}
natural in $M$, respectively called the \textbf{projection}, the \textbf{zero}, and the \textbf{sum}, which equip the tangent bundle with an additive structure on the fibre, the \textbf{vertical lift}, and finally, the \textbf{canonical flip}. The projection $p_M\colon\T M\to M$ is interpreted as the map that sends each tangent vector to its base point; the zero $z_M\colon M\to\T M$ corresponds to the map that sends each point $x$ of $M$ to the zero tangent vector $0_x$ at $x$; the sum $s_M\colon\T_2M\to\T M$ sums together two tangent vectors $u$ and $v$ that share the same base point $x$, where $\T_2M$ denotes the pullback
\begin{equation*}
\begin{tikzcd}
{\T_2M} & {\T M} \\
{\T M} & M
\arrow["{\pi_2}", from=1-1, to=1-2]
\arrow["{\pi_1}"', from=1-1, to=2-1]
\arrow["\lrcorner"{anchor=center, pos=0.125}, draw=none, from=1-1, to=2-2]
\arrow["{p_M}", from=1-2, to=2-2]
\arrow["{p_M}"', from=2-1, to=2-2]
\end{tikzcd}
\end{equation*}
and is a tangent pullback, in that it is preserved by all iterates $\T^n\=\T\T{\dots}\T$ of the tangent bundle functor. It is sometimes asked that $\T M$ be equipped with a \textbf{negation}, $n_M\colon\T M\to\T M$, also natural in $M$, which sends each tangent vector $u$ to $-u$. In that case, we say that the tangent category \textbf{has negatives}.

The vertical lift $l_M\colon\T M\to\T\T M$ encodes \emph{local linearity} of the tangent bundle by turning the following diagram
\begin{equation*}
\begin{tikzcd}
{\T_2M} & {\T\T M} \\
M & {\T M}
\arrow["{\xi_M}", from=1-1, to=1-2]
\arrow["{\pi_1p_M}"', from=1-1, to=2-1]
\arrow["{\T p_M}", from=1-2, to=2-2]
\arrow["{z_M}"', from=2-1, to=2-2]
\end{tikzcd}
\end{equation*}
into a pullback, where:
\begin{align*}
&\xi_M\=(z_{\T M}\times l_M)\T s_M\colon\T_2M\to\T\T M
\end{align*}
This universal condition makes the fibres $\T_xM\=p^{-1}(x)$ of the tangent bundle $\T M$ of $M$ into \textit{linear spaces}, in that the tangent bundle $\T\T_xM$ of $\T_xM$ is isomorphic to the Cartesian product $\T_xM\times\T_xM$. This linearity condition is precisely the condition enjoyed by Euclidean spaces, that are $\Real^n$, for which $\T\Real^n\cong\Real^n\times\Real^n$. This expresses the idea that an object in a tangent category is a \emph{locally linear} space and a morphism is a \emph{locally linear} function.

Lastly, the canonical flip, $c_M\colon\T\T M\to\T\T M$, is an isomorphism of the double tangent bundle that encodes the symmetry of partial derivatives, in that $\partial_i\partial_jf=\partial_j\partial_if$.

The archetypal example is the category of finitely-dimensional smooth manifolds with the usual notion of the tangent bundle functor. However, tangent categories capture geometric models well beyond smooth manifolds. In~\cite{cruttwell:algebraic-geometry}, the category of affine schemes was shown to carry a tangent structure which captures some key aspects of algebraic geometry. This construction was later extended in~\cite{ikonicoff:operadic-algebras-tangent-cats} to algebras over an arbitrary algebraic operad. Convenient manifolds and $C^\infty$-rings also form tangent categories as well as the subcategory of microlinear objects of a model of Synthetic Differential Geometry (SDG) and every Cartesian Differential Category (CDC)~\cite{cockett:tangent-cats}. More recently, in~\cite{blohmann:elastic-diffeological-spaces}, Blohmann showed that \emph{elastic} diffeological spaces also form a tangent category.
\end{content}

\begin{content}[Lie groups in tangent categories]
A \textbf{Lie group} (\cite[Definition~2.13]{lanfranchi:lie-correspondence-tangent-cats}) in a Cartesian tangent category, that is, a tangent category $(\X,\TT)$ with terminal object $\*$ and binary products $A\times B$ preserved by the tangent bundle functor, is a group object $G$ with unit map $\e\colon\*\to G$, multiplication $\m\colon G\times G\to G$, and inverse map $\i\colon G\to G$, of $(\X,\TT)$, which admits tangent space at the unit, in that the pullback
\begin{equation*}
\begin{tikzcd}
{\T_\e G} & {\T G} \\
{\*} & G
\arrow["{\iota_\e}", from=1-1, to=1-2]
\arrow[from=1-1, to=2-1]
\arrow["\lrcorner"{anchor=center, pos=0.125}, draw=none, from=1-1, to=2-2]
\arrow["{p_G}", from=1-2, to=2-2]
\arrow["\e"', from=2-1, to=2-2]
\end{tikzcd}
\end{equation*}
of the tangent bundle $p_G\colon\T G\to G$ along the unit map $\e$ exists and is preserved by all iterates $\T^n$ of the tangent bundle functor.
\end{content}

\begin{content}[Lie functors]
Every group object $G$ in a Cartesian tangent category whose tangent bundle $p_G$ has negation $n_G\colon\T G\to\T G$, comes with a functor
\begin{align*}
&\g(-)\colon\X\to\LieAlg^\op
\end{align*}
called the \textbf{Lie functor} of $G$, that sends every object $A$ of $\X$ to the Lie algebra $\g(A)$ of \textbf{$A$-Parametrized Left-Invariant Vector Fields} (\textbf{$A$-\PLIVFs}\!) of $G$ (\cite[Definition~3.18]{lanfranchi:lie-correspondence-tangent-cats}), which are maps $X\colon G\times A\to\T G$ that satisfy the following two conditions:
\begin{equation*}
\begin{tikzcd}
{G\times A} & {\T G} \\
& G
\arrow["X", from=1-1, to=1-2]
\arrow["{\pi_1}"', from=1-1, to=2-2]
\arrow["{p_G}", from=1-2, to=2-2]
\end{tikzcd}\qquad
\begin{tikzcd}
{G\times G\times A} & {G\times A} \\
{\T(G\times G)} & {\T G}
\arrow["{\m\times\id_A}", from=1-1, to=1-2]
\arrow["{z_G\times X}"', from=1-1, to=2-1]
\arrow["X", from=1-2, to=2-2]
\arrow["{\T\m}"', from=2-1, to=2-2]
\end{tikzcd}
\end{equation*}
The Lie functor $\g(-)$ of a Lie group is always representable (\cite[Theorem~4.20]{lanfranchi:lie-correspondence-tangent-cats}), and it is represented by the tangent space $\T_\e G$ of $G$. Thus, via the Yoneda lemma, the Lie algebra structure of each $\g(A)$ induces an \emph{internal} Lie algebra structure on $\T_\e G$ (\cite[Theorem~5.4]{lanfranchi:lie-correspondence-tangent-cats}).
\end{content}

\begin{content}[Differential Lie algebras]
The resulting Lie algebra object, denoted by $\g\in(\X,\TT)$, becomes a \textbf{differential Lie algebra} (\cite[Definition~5.1 and Theorem~5.12]{lanfranchi:lie-correspondence-tangent-cats}), in that the underlying object of $\g$ is a differential object, the additive structure of the Lie algebra coincides with that of the differential object $\g$, and the Lie bracket $[,]\colon\g\times\g\to\g$ is bilinear, that is:
\begin{equation*}
\begin{tikzcd}
{\T(\g\times\g)} & {\T(\g)} & {\T(\g\times\g)} \\
{\g\times\g} & \g & {\g\times\g}
\arrow["{\T[,]}", from=1-1, to=1-2]
\arrow["{\T[,]}"', from=1-3, to=1-2]
\arrow["{z_\g\times\lambda_\g}", from=2-1, to=1-1]
\arrow["{[,]}"', from=2-1, to=2-2]
\arrow["{\lambda_\g}"', from=2-2, to=1-2]
\arrow["{\lambda_\g\times z_\g}"', from=2-3, to=1-3]
\arrow["{[,]}", from=2-3, to=2-2]
\end{tikzcd}
\end{equation*}
where $\lambda_\g\colon\g\to\T\g$ is the vertical lift of the differential object $\g$.

Recall that a \textbf{differential object} (\cite[Definition~3.1]{cockett:differential-bundles}) in a Cartesian tangent category consists of a commutative monoid $(\A,0,+)$ equipped with a \textbf{vertical lift} $\lambda_\A\colon\A\to\T\A$, compatible with the additive structure, that makes the following morphism
\begin{align*}
&\hat\xi_\A\colon\A\times\A\xrightarrow{z_\A\times\lambda_\A}\T(\A\times\A)\xrightarrow{\T+}\T\A
\end{align*}
into an isomorphism. Equivalently, a differential object consists of a commutative monoid $(\A,0,+)$ equipped with a \textbf{differential projection}, which is a map $\hat p_\A\colon\T\A\to\A$ that is compatible with the additive structure of $\A$ and of its tangent bundle, and that turns the following diagram
\begin{equation*}
\begin{tikzcd}
\A & {\T\A} & \A
\arrow["{p_\A}"', from=1-2, to=1-1]
\arrow["{\hat p_\A}", from=1-2, to=1-3]
\end{tikzcd}
\end{equation*}
into a product diagram.
\end{content}

\begin{content}[Goal of the paper]
In this paper, we introduce and study general linear group objects and their Lie algebras in tangent categories. We believe this work is not only important in itself since it introduces an important class of group objects in the general context of tangent categories, but it is also preliminary to developing a rich theory of geometric constructions in the tangent-categorical world. In particular, general linear groups are an important ingredient in the definition of the \emph{frame} bundle of a vector bundle $\q\colon E\to M$, which is defined as a $\GL_n(\Real)$-principal bundle, where $n$ is the dimension of the typical fibre of $\q$.
\end{content}

\begin{content}[Structure of the paper]
The paper is divided into three major parts. In the first part, we take a deep dive into multilinearity in the context of tangent categories and, in particular, in the notion of internal multilinear Homs (Definition~\ref{definition:linearly-closed-tangent-category}), which provide internal Homs for multilinear maps between differential objects. While the notions of multilinearity are not new in the tangent-categorical literature (see~\cite{cruttwell:sector-forms-in-CDCs} and~\cite{lemay:jacobian}), so far, it has only been considered in the context of Cartesian differential categories. In this part, we will also introduce \textbf{Linear Endomorphism Classifiers} (Definition~\ref{definition:linear-endomorphism-classifier}), which are those internal multilinear Homs of type $\lnr[\A,\A]$. We will also study their algebraic properties and show that their maximal differential subring $\lnr_\A^-$ (Definition~\ref{definition:maximal-differential-subring}) exhibits the structure of a differential Lie algebra.

In Section~\ref{section:general-linear-groups}, we make use of the knowledge we gained in the first section about internal multilinear Homs and introduce general linear group objects, as classifiers of a special class of maps that we call \textbf{Kleisli parametrized linear automorphisms} (Definition~\ref{definition:kleisli-PLAs}), which are maps $f\colon\A\times M\to\T^n\A$ that lift to the Kleisli tangent category $\KL[\T_M^n,c_M^n]$ as linear automorphisms $f\colon\A\nto\A$ (Lemma~\ref{lemma:kleisli-PLI}). After introducing general linear groups $\GL_\A$ of each differential object $\A$, we show how to construct them directly from the linear endomorphism classifier $\lnr_\A$ via a tangent triple equalizer (Theorem~\ref{theorem:triple-equalizer}).

In the third and last part, we study the differential Lie algebras of the general linear groups. We begin by constructing a natural Lie algebra homomorphism $\chi_M\colon\gl_\A(M)\to\lnr_\A^-(M)$ between the Lie functor of $\GL_\A$ and the functor $\lnr_\A^-(-)$, which is represented by the Lie algebra $\lnr_\A^-$ of the maximal subring of the linear endomorphism classifier (Proposition~\ref{proposition:the-map-chi}). Using this construction and the characterization of Theorem~\ref{theorem:triple-equalizer}, we prove the main result of the paper: that the differential Lie algebra $\gl_\A$ of the general linear group $\GL_\A$ of a differential object $\A$ is isomorphic to the differential Lie algebra $\lnr_\A^-$ (Theorem~\ref{theorem:main-theorem}). Finally, when $\A$ admits negatives, this theorem recaptures the classic result of differential and algebraic geometry, that $\gl_\A\cong\lnr_\A$ (Corollary~\ref{corollary:main}).
\end{content}

\begin{conventions}
We shall denote composition diagrammatically; that is, $fg$ is the composition of $f$ followed by $g$. When considering functors or concrete functions, we switch to the functional notation, so for example, we shall denote by $GFM$ or $GFf$ the \emph{application} of the functor $G$ to the application of the functor $F$ to the object $M$ and to the morphism $f$.

We shall denote generic differential objects as $\A,\B,\C$ and we shall denote the zero, sum, vertical lift, and differential projection of $\A$ by $0_\A$, $+_\A$, $\lambda_\A$, and $\hat p_\A$. We shall omit the subscript $\A$ when it will be clear from the context.
\end{conventions}


\section{Multilinearity in tangent categories}
\label{section:linearity}
To introduce general linear group objects in tangent categories, we first take a detour and study internal multilinear Homs.

\subsection{A fibrational perspective on multilinearity}
\label{subsection:fibrational-multilinearity}
In~\cite{lemay:jacobian}, Lemay introduced \emph{linearly closed Cartesian differential categories} to give a synthetic perspective on the Jacobian matrix of smooth morphisms. In this section, we extend this concept to Cartesian tangent categories. We begin by adapting to our context the notions of linear and multilinear morphisms, introduced by Cruttwell in~\cite[Definition~2.9]{cruttwell:sector-forms-in-CDCs} in the context of Cartesian differential categories.

\begin{definition}[{\cite[Definition~2.9, adapted for tangent categories]{cruttwell:sector-forms-in-CDCs}}]
\label{definition:multilinearity}
In a Cartesian tangent category:
\begin{enumerate}
\item A morphism $f\colon M_1\ntimes\A_k\ntimes M_n\to\B$ is \textbf{linear in the $k$-th argument}, where $\A_k$ and $\B$ are differential objects, provided that
\begin{equation*}
\begin{tikzcd}
{\T(M_1\ntimes\A_k\ntimes M_n)} & {\T\B} \\
{M_1\ntimes\A_k\ntimes M_n} & \B
\arrow["{\T f}", from=1-1, to=1-2]
\arrow["{\lambda_k}", from=2-1, to=1-1]
\arrow["f"', from=2-1, to=2-2]
\arrow["{\lambda}"', from=2-2, to=1-2]
\end{tikzcd}
\end{equation*}
where $\lambda_k\=z_{M_1}\ntimes\lambda_{\A_k}\ntimes z_{M_n}$.

\item A morphism $f\colon\A_1\ntimes\A_n\to\B$ is $(\A_1\,\A_n)$-\textbf{multilinear} (or simply \textbf{multilinear}), where $\A_1\,\A_n,\B$ are differential objects, provided that $f$ is linear in each $k$-th argument, for $k=1\,n$.

\item An $(\A)$-multilinear morphism $f\colon\A\to\B$ is called \textbf{linear}\footnote{This corresponds to the usual definition of linearity between differential objects introduced in~\cite{cockett:differential-bundles}.}.

\item An $(\A_1,\A_2)$-multilinear morphism $f\colon\A_1\times\A_2\to\B$ is called \textbf{bilinear}.
\end{enumerate}
\end{definition}

Differential objects together with multilinear maps form a multicategory denoted by $\DO_\multi[\X,\TT]$.

We shall spend some time reinterpreting multilinearity through a fibrational lens. Our goal is to express multilinearity as proper \emph{linearity at each entry}. We begin by recalling that the \textbf{simple fibration} of a Cartesian category $\X$ (introduced in~\cite[Section~1.3]{jacobs:fibrations-categorical-logic-type-theory}) is the (cloven) fibration
\begin{align*}
&\Pi\colon\Simple[\X]\to\X
\end{align*}
whose total category $\Simple[\X]$ is the category of pairs $(M,A)$ of objects of $\X$, whose morphisms are pairs $(f,\varphi)$ of morphisms $f\colon M\to N$ and $\varphi\colon A\times M\to B$\footnote{We adopted the non-standard convention of switching between $A$ and $M$. This convention is more natural when we consider left-invariant vector fields.}. The identity morphism of $(M,A)$ is given by $(\id_M,\pi_1)$ and composition is defined by the following formula:
\begin{align*}
&(f,\varphi)(g,\psi)=\left(fg,A\times M\xrightarrow{\id_A\times\Delta_M}A\times M\times M\xrightarrow{\varphi\times f}B\times N\xrightarrow{\psi}C\right)
\end{align*}
The fibration $\Pi$ sends each $(M,A)$ to $M$ and each $(f,\varphi)$ to $f$. When the Cartesian category $\X$ is equipped with a Cartesian tangent structure $\TT$, $\Pi$ becomes a (cloven) \textbf{tangent fibration}~\cite[Definition~5.2]{cockett:differential-bundles}. Recall that a (cloven) tangent fibration consists of a fibration between two tangent categories whose underlying functor preserves the tangent structure strictly and the tangent bundle functors preserve the (cloven) Cartesian lifts.

\begin{proposition}[{\cite[Proposition~3.7]{lanfranchi:lie-correspondence-tangent-cats}}]
\label{proposition:simple-tangent-fibration}
The simple fibration $\Pi\colon\Simple[\X,\TT]\to(\X,\TT)$ of a Cartesian tangent category is a cloven tangent fibration.
\end{proposition}

By~\cite[Theorem~5.3]{cockett:differential-bundles}, the fibres $\Simple_M[\X,\TT]$ become tangent categories and the substitution functors between the fibres become strong tangent morphisms. Concretely, the tangent bundle functor of $\Simple[\X,\TT]$ sends each $(M,A)$ to $(\T M,\T A)$ and each $(f,\varphi)$ to $(\T f,\varphi_\T)$, where:
\begin{align*}
&\varphi_\T\colon\T A\times\T M\xrightarrow{\cong}\T(A\times M)\xrightarrow{\T\varphi}\T B
\end{align*}
Furthermore, for each $M\in(\X,\TT)$, the fibre $\Simple_M[\X,\TT]$ over $M$ of the simple tangent fibration becomes a tangent category whose objects are the objects of $\X$, morphisms are morphisms $\varphi\colon A\times M\to B$, and the tangent bundle functor sends each $A$ to $\T A$ and each $\varphi$ to:
\begin{align*}
&\T A\times M\xrightarrow{\id_{\T A}\times z_M}\T(A\times M)\xrightarrow{\T\varphi}\T B
\end{align*}

\begin{lemma}[{\cite[Lemma~4.13]{lanfranchi:lie-correspondence-tangent-cats}}]
\label{lemma:inclusion-functor-preserves-differential-objects}
The inclusion functor $(\X,\TT)\to\Simple_M[\X,\TT]$ that sends each $A$ to itself and each $f\colon A\to B$ to $f^\KL$ defined by $\SQ{f^\KL}\=\pi_1f\colon A\times M\to B$ becomes a strict tangent morphism that preserves tangent limits. Moreover, it preserves differential objects.
\end{lemma}

We shall denote by $\A_1\ntimes\hat\A_k\ntimes\A_n$ the Cartesian product of a list of differential objects in which the $k$-th term has been omitted.

\begin{lemma}
\label{lemma:multilinearity-fibrationally}
A morphism $f\colon M_1\ntimes\A_k\ntimes M_n\to\B$ is linear in the $k$-th argument if and only if
\begin{align*}
\SQ{f_k}\colon\A_k\times M_1\ntimes\hat\A_k\ntimes M_n\xrightarrow{\Shuffle_k}M_1\ntimes\A_k\ntimes M_n\xrightarrow{f}\B
\end{align*}
lifts to a linear morphism $\A_k\nto\B$ in $\Simple_{M_1\ntimes\hat\A_k\ntimes M_n}[\X,\TT]$, where $\Shuffle_k$ denotes the permutation that sends the first entry to the $k$-th entry and preserves the orders of the other entries.
\end{lemma}
\begin{proof}
Consider a morphism $f\colon M_1\ntimes\A_k\ntimes M_n\to\B$ and suppose that $f_k\colon\A_k\nto\B$ is linear in $\Simple_k[\X,\TT]\=\Simple_{M_1\ntimes\hat\A_k\ntimes M_n}[\X,\TT]$. By definition, this means that $f_k\lambda_\B=\lambda_{\A_k}\T_{(k)}f_k$ in $\Simple_k[\X,\TT]$. By construction, the tangent bundle functor $\T_{(k)}$ sends $f_k$ to:
\begin{align*}
&\SQ{\T_{(k)}f_k}=(\id\times z\ntimes z)\T f_k=(\id\times z\ntimes z)\Shuffle_k\T f=\\
&\qquad=\Shuffle_k(z\ntimes\id\ntimes z)\T f
\end{align*}
Thus, using the formula that defines the coKleisli composition in $\Simple_k[\X,\TT]$, we rewrite the linearity condition as follows:
\begin{align*}
&\Shuffle_k f\lambda_\B=\Shuffle_k(z\ntimes\lambda_{\A_k}\ntimes z)\T f=\Shuffle_k\lambda_k\T f
\end{align*}
However, since $\Shuffle_k$ is just a permutation, it is invertible, thus $fl\lambda_\B=\lambda_k\T f$, which means that $f$ is linear in the $k$-th argument.
\end{proof}

As explained in the introduction, a differential object is a commutative monoid $\A$ equipped with a vertical lift $\lambda\colon\A\to\T\A$ that is used to construct an isomorphism $\xi\colon\A\times\A\to\T\A$. Equivalently, a differential object is a commutative monoid $\A$ equipped with a differential projection $\hat p$ that turns $\T\A$ into the product $\A\times\A$. A morphism $f\colon\A\to\B$ is linear if $f\lambda_\B=\lambda_\A\T f$, or as expressed via the differential projections, if $\T f\hat p_\B=\hat p_\A f$. Thanks to Lemma~\ref{lemma:multilinearity-fibrationally}, we can rephrase multilinearity as a compatibility condition with the differential projections.

\begin{lemma}
\label{lemma:multilinearity-as-p-hat}
A morphism $f\colon M_1\ntimes\A_k\ntimes M_n\to\B$ is linear in its $k$-th argument if and only if the following diagram
\begin{equation*}
\begin{tikzcd}
{M_1\ntimes\T\A_k\ntimes M_n} && {\T(M_1\ntimes\A_k\ntimes M_n)} & {\T\B} \\
{M_1\ntimes\A_k\ntimes M_n} &&& \B
\arrow["{z\ntimes z\times\id\times z\ntimes z}", from=1-1, to=1-3]
\arrow["{\hat p_k}"', from=1-1, to=2-1]
\arrow["{\T f}", from=1-3, to=1-4]
\arrow["{\hat p_\B}", from=1-4, to=2-4]
\arrow["f"', from=2-1, to=2-4]
\end{tikzcd}
\end{equation*}
commutes, where:
\begin{align*}
&\hat p_k\=\id\ntimes\id\times\hat p_{\A_k}\times\id\ntimes\id
\end{align*}
In particular, a morphism $f\colon\A_1\ntimes\A_n\to\B$ is multilinear if and only if, for each $k=1\,n$ the following diagram
\begin{equation*}
\begin{tikzcd}
{\A_1\ntimes\T\A_k\ntimes \A_n} && {\T(\A_1\ntimes\A_n)} & {\T\B} \\
{\A_1\ntimes\A_n} &&& \B
\arrow["{z\ntimes z\times\id\times z\ntimes z}", from=1-1, to=1-3]
\arrow["{\hat p_k}"', from=1-1, to=2-1]
\arrow["{\T f}", from=1-3, to=1-4]
\arrow["{\hat p_\B}", from=1-4, to=2-4]
\arrow["f"', from=2-1, to=2-4]
\end{tikzcd}
\end{equation*}
commutes.
\end{lemma}
\begin{proof}
A morphism $f\colon\A\to\B$ is linear if and only if it satisfies $\T f\hat p_\B=\hat p_\A f$. From Lemma~\ref{lemma:multilinearity-fibrationally}, a morphism $f\colon M_1\ntimes\A_k\ntimes M_n\to\B$ is linear in its $k$-th entry if and only if $f_k\colon\A_k\nto\B$, defined by $\SQ{f_k}\=\Shuffle_kf$, becomes a linear morphism in $\Simple_k[\X,\TT]\=\Simple_{M_1\ntimes\hat\A_k\ntimes M_n}[\X,\TT]$. However, $f_k$ is linear if and only if $\T_{(k)}f_k\hat p_\B=\hat p_{\A_k}f_k$ in $\Simple_k[\X,\TT]$, where $\T_{(k)}$ denotes its tangent bundle functor. However, by unwrapping the definition of $\T_{(k)}$ and of the composition, we rewrite it as follows:
\begin{align*}
&\SQ{\T_{(k)}f_k\hat p_\B}=(\id\times z\ntimes z)\T f_k\hat p_\B        &&\SQ{\hat p_\A f_k}=(\hat p_{\A_k}\times\id\ntimes\id)f_k
\end{align*}
Therefore, by precomposing by $\Shuffle_k^{-1}$, we obtain the desired equation:
\begin{align*}
&(z\ntimes\id\ntimes z)\T f\hat p_\B=(\id\ntimes\hat p_{\A_k}\ntimes\id)f=\hat p_kf\qedhere
\end{align*}
\end{proof}

In the next section, we will extensively consider morphisms $f\colon\A_1\ntimes\A_n\times M\to\B$ which are linear in each $k$-th entry, for $k=1\,n$. Thanks to the fibrational perspective we presented, we can reinterpret such maps as \emph{multilinear} morphisms $f\colon\A_1\ntimes\A_n\nto\B$ in the tangent category $\Simple_M[\X,\TT]$.

\begin{definition}
\label{definition:M-parametrized-multilinear-map}
An \textbf{$M$-parametrized multilinear map} is a map $f\colon\A_1\ntimes\A_n\times M\to\B$ that lifts to a multilinear morphism $f\colon\A_1\,\A_n\to\B$ in $\Simple_M[\X,\TT]$.
\end{definition}

\subsection{Multilinearly closed tangent categories}
\label{subsection:internal-linear-homs}
We shall now introduce internal multilinear Homs. We begin by adapting Lemay's definition of internal \emph{linear} Homs in the context of tangent categories and extend it to capture multilinear maps.

\begin{definition}[{\cite[Definition~3.1, adapted for tangent categories]{lemay:jacobian}}]
\label{definition:linearly-closed-tangent-category}
In a Cartesian tangent category, given a list of differential objects $\A_1\,\A_n,\B$, the \textbf{internal multilinear Hom} from $\A_1\,\A_n$ to $\B$ consists of an object $\lnr[\A_1\,\A_n;\B]$ equipped with a morphism
\begin{align*}
&\ev\colon\A_1\ntimes\A_n\times\lnr[\A_1\,\A_n;\B]\to\B
\end{align*}
called the \textbf{evaluation map}, subject to the following conditions:
\begin{description}
\item[ILH.1] $\lnr[\A_1\,\A_n;\B]$ is a differential object and the evaluation map is multilinear, in that:
\begin{equation*}
\begin{tikzcd}
{\T(\A_1\ntimes\A_n\times\lnr[\A_1\,\A_n;\B])} & {\T\B} \\
{\A_1\ntimes\A_n\times\lnr[\A_1\,\A_n;\B]} & \B
\arrow["{\T\ev}", from=1-1, to=1-2]
\arrow["{\lambda_k}", from=2-1, to=1-1]
\arrow["{\ev}"', from=2-1, to=2-2]
\arrow["{\lambda_\B}"', from=2-2, to=1-2]
\end{tikzcd}
\end{equation*}
for $k=1\,n+1$;

\item[ILH.2] The evaluation map classifies multilinear maps, in that for every $M$-parametrized multilinear map $f\colon\A_1\ntimes\A_n\times M\to\B$, there exists a unique morphism $\hat f\colon M\to\lnr[\A_1\,\A_n;\B]$ that makes the following diagram
\begin{equation*}
\begin{tikzcd}
{\A_1\ntimes\A_n\times\lnr[\A_1\,\A_n;\B]} & \B \\
{\A_1\ntimes\A_n\times M}
\arrow["{\ev}", from=1-1, to=1-2]
\arrow["{\id\times\hat f}", from=2-1, to=1-1]
\arrow["f"', from=2-1, to=1-2]
\end{tikzcd}
\end{equation*}
commutative.
\end{description}
When $n=1$, we refer to $\lnr[\A;\B]$ as the \textbf{internal linear Hom}. A tangent category is \textbf{multilinearly closed} provided that the internal multilinear Hom $\lnr[\A_1\,\A_n;\B]$ exists for each list $\A_1\,\A_n,\B$ of differential objects.
\end{definition}

\begin{remark}
\label{remark:internal-linear-homs}
Our Definition~\ref{definition:linearly-closed-tangent-category} and Lemay's original \cite[Definition~3.1]{lemay:jacobian} differ in a few details. First, the ambient category is not simply a Cartesian differential category but, in fact, a Cartesian tangent category. Second, the universal property \textbf{[ILH.2]}. In fact, this universal property classifies not only multilinear maps but also every \emph{parametrized} multilinear map, where the parameter object $M$ can be chosen arbitrarily in the tangent category. Finally, the original definition only defines internal \emph{linear} Homs, while our version extends this to general multilinear maps.
\end{remark}

Despite these differences between ours and Lemay's original definition, they are equivalent in the category of differential objects. We shall omit the easy proof of this statement.

\begin{proposition}
\label{proposition:internal-linear-hom-js}
The Cartesian differential category $\DO[\X,\TT]$ of differential objects of a multilinearly closed tangent category $(\X,\TT)$ is linearly closed in the sense of~\cite[Definition~3.1]{lemay:jacobian}. In particular, each $\lnr[\A,\B]$ becomes an internal linear Hom as in~\cite[Definition~3.1]{lemay:jacobian}.
\end{proposition}

We now show that each $\lnr[\A_1\,\A_n;\B]$ is an exponential object in the multicategory $\DO_\multi[\X,\TT]$ of differential objects and multilinear morphisms. Recall that an exponential object in a multicategory consists of an object $\B^{\A_1\,\A_n}$ equipped with a multimorphism $\ev\colon\A_1\,\A_n,\B^{\A_1\,\A_n}\to\B$ subject to the following universal property. For every multimorphism $f\colon\A_1\,\A_n,\A'_1\,\A'_m\to\B$, there exists a unique morphism $\hat f\colon\A_1'\,\A_m'\to\B^{\A_1\,\A_n}$ such that $(\id_{\A_1}\,\id_{\A_n},\hat f)\ev=f$.

\begin{proposition}
\label{proposition:exponential-objects}
Internal multilinear Homs are exponential objects in the multicategory $\DO_\multi[\X,\TT]$.
\end{proposition}
\begin{proof}
By \textbf{[ILH.1]}, each internal multilinear Hom is an object of $\DO_\multi[\X,\TT]$ and the evaluation map $\ev$ becomes a multimorphism $\ev\colon\A_1\,\A_n,\lnr[\A_1\,\A_n;\B]\to\B$ of $\DO_\multi[\X,\TT]$. It is left only to prove that for each $f\colon\A_1\,\A_n,\A_1'\,\A_m'\to\B$, $\hat f\colon\A_1'\ntimes\A_m'\to\lnr[\A_1\,\A_n;\B]$, which is the unique morphism satisfying $(\id\times\hat f)\ev=f$, is multilinear. To this end, we shall compute:
\begin{align*}
(\id\ntimes\id\times(\hat p_k\hat f))\ev&=~(\id\ntimes\id\times\hat p_k)(\id\ntimes\id\times\hat f)\ev\\
&=~(\id\ntimes\id\times\hat p_k)f                   \tag{Definition of $\hat f$}\\
&=~(z\ntimes z\times\id)\T f\hat p                  \tag{$f$ is multilinear}\\
&=~(z\ntimes z\times\id)(\id\ntimes\id\times\T\hat f)\T\ev\hat p    \tag{Definition of $\hat f$}\\
&=~(\id\ntimes\id\times\T\hat f)(z\ntimes z\times\id)\T\ev\hat p\\
&=~(\id\ntimes\id\times\T\hat f)(\id\ntimes\id\times\hat p)\ev          \tag{$\ev$ is multilinear}\\
&=~(\id\ntimes\id\times(\T\hat f\hat p))\ev
\end{align*}
By invoking \textbf{[ILH.2]}, we conclude that $\hat p_k\hat f=\T\hat f\hat p$, proving, thanks to Lemma~\ref{lemma:multilinearity-as-p-hat}, that $\hat f$ is in fact multilinear.
\end{proof}

In Definition~\ref{definition:linearly-closed-tangent-category}, \textbf{[ILH.1]} establishes that the object $\lnr[\A_1\,\A_n;\B]$ should be a differential object and that the evaluation map should be multilinear, while \textbf{[ILH.2]} establishes the class of morphisms that internal multilinear Homs classify. It turns out that one can \emph{derive} that $\lnr[\A_1\,\A_n;\B]$ is a differential object by enlarging the class of morphisms classified by the evaluation map.

\begin{theorem}
\label{theorem:internal-linear-homs-classify-kleisli-multilinear-morphisms}
Given a list $\A_1\,\A_n,\B$ of differential object, consider an object $\lnr[\A_1\,\A_n;\B]$ equipped with a map $\ev\colon\A_1\ntimes\A_n\times\lnr[\A_1\,\A_n;\B]\to\B$. Then $\lnr[\A_1\,\A_n;\B]$ is an internal multilinear Hom with evaluation map $\ev$ if and only if the following conditions are satisfied:
\begin{description}
\item[ILH.1'] The map $\ev$ is a $\lnr[\A_1\,\A_n;\B]$-parametrized multilinear map;

\item[ILH.2'] For every $M$-parametrized multilinear map $f\colon\A_1\ntimes\A_n\times M\to\T\B$, there exists a unique morphism $\hat f\colon M\to\T\lnr[\A_1\,\A_n;\B]$ that makes the following diagram
\begin{equation*}
\begin{tikzcd}
{\T(\A_1\ntimes\A_n\times\lnr[\A_1\,\A_n;\B])} & {\T\B} \\
{\A_1\ntimes\A_n\times M}
\arrow["{\T\ev}", from=1-1, to=1-2]
\arrow["{z\ntimes z\times\hat f}", from=2-1, to=1-1]
\arrow["f"', from=2-1, to=1-2]
\end{tikzcd}
\end{equation*}
commutative.
\end{description}
\end{theorem}
\begin{proof}
We begin by assuming that $\lnr[\A_1\,\A_n;\B]$ is an internal multilinear Hom object and we prove that it satisfies \textbf{[ILH.1']} and \textbf{[ILH.2']}. \textbf{[ILH.1']} follows immediately from \textbf{[ILH.1]}. To prove \textbf{[ILH.2']}, consider a parametrized multilinear map $f\colon\A_1\ntimes\A_n\times M\to\T\B$. Since $\B$ is a differential object, $\T\B\cong\B\times\B$, so we can consider $f_1\=fp_\B$ and $f_2\=f\hat p_\B$, which are also parametrized multilinear maps. Therefore, by \textbf{[ILH.2]}, we find two morphisms $\hat f_1,\hat f_2\colon M\to\lnr[\A_1\,\A_n;\B]$, satisfying:
\begin{align*}
&(\id\times\hat f_1)\ev=f_1             &&(\id\times\hat f_2)\ev=f_2
\end{align*}
However, by \textbf{[ILH.1]}, $\T\lnr[\A_1\,\A_n;\B]\cong\lnr[\A_1\,\A_n;\B]\times\lnr[\A_1\,\A_n;\B]$. So, we can define $\hat f\colon M\to\T\lnr[\A_1\,\A_n;\B]$ as the morphism $\hat f=\<\hat f_1,\hat f_2\>$. We shall prove that $(z\ntimes z\times\hat f)\T\ev=f$:
\begin{align*}
(z\ntimes z\times\hat f)\T\ev p&=~(zp\ntimes zp\times\hat fp)\ev            \tag{Naturality of $p$}\\
&=~(\id\times\hat f_1)\ev                                                   \tag{$zp=\id$, $\hat fp=\hat f_1$}\\
&=~f_1                                                                      \tag{Definition of $\hat f_1$}
\end{align*}
Moreover:
\begin{align*}
(z\ntimes z\times\hat f)\T\ev\hat p&=~(\id\times\hat f\hat p)\ev            \tag{$\ev$ is multilinear}\\
&=~(\id\times\hat f_2)\ev                                                   \tag{$\hat f\hat p=\hat f_2$}\\
&=~f_2                                                                      \tag{Definition of $\hat f_2$}
\end{align*}
This proves that $(z\ntimes z\times\hat f)\T\ev=f$. Now, consider a morphism $g\colon M\to\T\lnr[\A_1\,\A_n;\B]$ such that $(z\ntimes z\times g)\T\ev=f$. Thus:
\begin{align*}
&(\id\times\hat f_1)\ev=f_1=fp=(z\ntimes z\times g)\T\ev p=(\id\times gp)\ev\\
&(\id\times\hat f_2)\ev=f_2=f\hat p=(z\ntimes z\times g)\T\ev\hat p=(\id\times g\hat p)\ev
\end{align*}
Therefore, by \textbf{[ILH.2]}, $\hat f_1=gp$ and $\hat f_2=g\hat p$; thus, $g=\hat f$. This proves \textbf{[ILH.2']}.

Conversely, we now assume that $\lnr[\A_1\,\A_n;\B]$ with $\ev$ satisfies \textbf{[ILH.1']} and \textbf{[ILH.2']} and we shall prove that $\lnr[\A_1\,\A_n;\B]$ is an internal multilinear Hom object. We begin by proving that \textbf{[ILH.2]} holds. To this end, consider an $M$-parametrized multilinear map $f\colon\A_1\ntimes\A_n\times M\to\B$. It is easy to see that $fz$ is also parametrized multilinear. Thus, by \textbf{[ILH.2']}, there exists a unique morphism $g\colon M\to\T\lnr[\A_1\,\A_n;\B]$ such that $(z\ntimes z\times g)\T\ev=fz$. We shall prove that $gp$ satisfies $(\id\times gp)\ev=f$:
\begin{align*}
f&=~fzp                                         \tag{$zp=\id$}\\
&=~(z\ntimes z\times g)\T\ev p                  \tag{Definition of $g$}\\
&=~(\id\times gp)\ev                            \tag{Naturality of $p$}
\end{align*}        
Now, consider another morphism $h\colon M\to\lnr[\A_1\,\A_n;\B]$ that satisfies $(\id\times h)\ev=f$. Thus:
\begin{align*}
&(z\ntimes z\times g)\T\ev=fz=(\id\times h)\ev z=(z\ntimes z\times hz)\T\ev
\end{align*}
where we used the naturality of $z$. Thus, by \textbf{[ILH.2']}, $hz=g$; $h=gp$, proving uniqueness. This shows that \textbf{[ILH.2]} holds. Thus, in the following, we shall assume that \textbf{[ILH.2]} is given. Next, we show that $\lnr[\A_1\,\A_n;\B]$ is a differential object and that $\ev$ is multilinear. We begin by constructing the additive structure of $\lnr[\A_1\,\A_n;\B]$. We define the zero $0\colon\*\to\lnr[\A_1\,\A_n;\B]$ as the unique morphism that satisfies the following equation: $(\id\times0)\ev=!0$. It is not hard to see that $!0$ is a parametrized multilinear map. Thus, the existence of $0$ follows from the universal property \textbf{[ILH.2]}. Similarly, we define the sum $+\colon\lnr[\A_1\,\A_n;\B]\times\lnr[\A_1\,\A_n;\B]\to\lnr[\A_1\,\A_n;\B]$ as the unique morphism that makes the following diagram
\begin{equation*}
\adjustbox{width=\linewidth}{
\begin{tikzcd}
{\A\times\lnr[\A_1\,\A_n;\B]} && \B \\
&& {\B\times\B} \\
{\A_1\ntimes\A_n\times\lnr[\A_1\,\A_n;\B]\times\lnr[\A_1\,\A_n;\B]} && {\A_1\ntimes\A_n\times\lnr[\A_1\,\A_n;\B]\times\A_1\ntimes\A_n\times\lnr[\A_1\,\A_n;\B]}
\arrow["\ev", from=1-1, to=1-3]
\arrow["{+}"', from=2-3, to=1-3]
\arrow["{\id\times+}", from=3-1, to=1-1]
\arrow["{(\Delta\times\id)\ex}"', from=3-1, to=3-3]
\arrow["{\ev\times\ev}"', from=3-3, to=2-3]
\end{tikzcd}
}
\end{equation*}
commutative, where $\Delta\colon M\to M\times M$ is the diagonal map and $\ex\colon M\times M\times N\times N\to M\times N\times M\times N$ denotes the  map that exchanges the internal terms. To check that this morphism is well-defined, we shall prove that $(\Delta\times\id)\ex(\ev\times\ev)+$ is a parametrized multilinear map by computing the following:
\begin{align*}
(\Delta\times\id)\ex(\ev&\times\ev)+\lambda\\
&=~(\Delta\times\id)\ex((\ev\lambda)\times(\ev\lambda))\T+     \tag{Additivity of $\lambda$}\\
&=~\<\lambda_k,z,\lambda_k,z\>\T((\ev\times\ev)+)      \tag{\textbf{[ILH.1']}}\\
&=~(\lambda_k\times z)\T((\Delta\times\id)\ex(\ev\times\ev)+)
\end{align*}
We leave it to the reader to prove that $(\lnr[\A_1\,\A_n;\B],0,+)$ forms a commutative monoid, which in the internal language is equivalent to proving the following equations:
\begin{align*}
&(f(x_1\,x_n)+g(x_1\,x_n))+h(x_1\,x_n)=f(x_1\,x_n)+(g(x_1\,x_n)+h(x_1\,x_n))\\
&f(x_1\,x_n)+g(x_1\,x_n)=g(x_1\,x_n)+f(x_1\,x_n)\\
&f(x_1\,x_n)+0=f(x_1\,x_n)
\end{align*}
where $f,g$, and $h$ are \emph{internal morphisms} in $\lnr[\A_1\,\A_n;\B]$. We shall define the differential projection $\hat p$ as the unique morphism that makes the following diagram
\begin{equation*}
\adjustbox{width=\linewidth}{
\begin{tikzcd}
{\A_1\ntimes\A_n\times\T\lnr[\A_1\,\A_n;\B]} & {\T(\A_1\ntimes\A_n\times\lnr[\A_1\,\A_n;\B])} & {\T\B} \\
{\A_1\ntimes\A_n\times\lnr[\A_1\,\A_n;\B]} && \B
\arrow["{z\times\id}", from=1-1, to=1-2]
\arrow["{\id\times\hat p}"', from=1-1, to=2-1]
\arrow["{\T\ev}", from=1-2, to=1-3]
\arrow["{\hat p}", from=1-3, to=2-3]
\arrow["{\ev}"', from=2-1, to=2-3]
\end{tikzcd}
}
\end{equation*}
commutative. Since $z$ and $\hat p$ are linear, by composition, $(z\times\id)\T\ev\hat p$ is a parametrized multilinear map. Thus, $\hat p$ is well-defined. We leave it to the reader to prove the compatibilities between $\hat p$ and the additive structure of $\lnr[\A_1\,\A_n;\B]$, which is a consequence of the compatibilities between the differential projection $\hat p$ and the additive structure of $\B$. Instead, we focus on proving the universal property that establishes that $\T\lnr[\A_1\,\A_n;\B]$ is isomorphic to $\lnr[\A_1\,\A_n;\B]\times\lnr[\A_1\,\A_n;\B]$ with projections $p$ and $\hat p$. Consider two morphisms $f,g\colon M\to\lnr[\A_1\,\A_n;\B]$. Our goal is to construct a unique morphism $\<f,g\>\colon M\to\T\lnr[\A_1\,\A_n;\B]$ that satisfies the following two equations: $\<f,g\>p=f$ and $\<f,g\>\hat p=g$. We begin by considering the following two morphisms:
\begin{align*}
&\tilde f\colon\A_1\ntimes\A_n\times M\xrightarrow{\id\times f}\A_1\ntimes\A_n\times\lnr[\A_1\,\A_n;\B]\xrightarrow{\ev}\B\\
&\tilde g\colon\A_1\ntimes\A_n\times M\xrightarrow{\id\times g}\A_1\ntimes\A_n\times\lnr[\A_1\,\A_n;\B]\xrightarrow{\ev}\B
\end{align*}
Since $\ev$ is a parametrized multilinear map, it follows that $\tilde f$ and $\tilde g$ are also parametrized multilinear maps. Using that $\B$ is a differential object, we can define a morphism $\<\tilde f,\tilde g\>\colon\A_1\ntimes\A_n\times M\to\T\B$ as the unique morphism that satisfies the following two equations: $\<\tilde f,\tilde g\>p=\tilde f$ and $\<\tilde f,\tilde g\>\hat p=\tilde g$. We shall prove that $\<\tilde f,\tilde g\>$ is a parametrized multilinear map by computing the following:
\begin{align*}
\<\tilde f,\tilde g\>\T\lambda c=\<\tilde f\lambda,\tilde g\lambda\>=(\lambda_k\times z_M)\T(\<\tilde f,\tilde g\>)
\end{align*}
where we used that $\T\lambda c$ becomes $\lambda\times\lambda$ when $\T\B$ is seen as a Cartesian product. Thus, since $\<\tilde f,\tilde g\>$ is a parametrized multilinear map, by \textbf{[ILH.2']}, there exists a unique $\<f,g\>\colon M\to\T\lnr[\A_1\,\A_n;\B]$ that satisfies the following equation:
\begin{align*}
&(z\times\<f,g\>)\T\ev=\<\tilde f,\tilde g\>
\end{align*}
We now prove that $\<f,g\>$ satisfes the desired equations:
\begin{align*}
(\id\times(\<f,g\>p))\ev&=~(z\times\<f,g\>)(p\times p)\ev  \tag{$zp=\id$}\\
&=~(z\times\<f,g\>)\T\ev p     \tag{Naturality of $p$}\\
&=~\<\tilde f,\tilde g\>p                    \tag{Definition of $\<f,g\>$}\\
&=~\tilde f                                   \tag{Definition of $\<\tilde f,\tilde g\>$}\\
&=~(\id\times f)\ev              \tag{Definition of $\tilde f$}
\end{align*}
Thus, by the universal property \textbf{[ILH.2]}, $\<f,g\>p=f$. Now, let us show that $\<f,g\>\hat p=g$:
\begin{align*}
(\id\times(\<f,g\>\hat p))\ev&=~(\id\times\<f,g\>)(\id\times\hat p)\ev\\
&=~(\id\times\<f,g\>)(z\times\id)\T\ev\hat p     \tag{Definition of $\hat p$}\\
&=~(z\times\<f,g\>)\T\ev\hat p\\
&=~\<\tilde f,\tilde g\>\hat p        \tag{Definition of $\<f,g\>$}\\
&=~\tilde g                                   \tag{Definition of $\<\tilde f,\tilde g\>$}\\
&=~(\id\times g)\ev              \tag{Definition of $\tilde g$}
\end{align*}
Thus, using again \textbf{[ILH.2]}, we conclude that $\<f,g\>\hat p=g$. Now, let $h\colon M\to\T\lnr[\A_1\,\A_n;\B]$ satisfy $hp=f$ and $h\hat p=g$. Thus, following the same steps as before, we can easily prove that:
\begin{align*}
&(z\times h)\T\ev p=(\id\times f)\ev=\tilde f    &&(z\times h)\T\ev\hat p=(\id\times g)\ev=\tilde g
\end{align*}
However, there is only one morphism that satisfies these two equations, that is:
\begin{align*}
&(z\times h)\T\ev=\<\tilde f,\tilde g\>=(\<z\times\<f,g\>)\T\ev
\end{align*}
By \textbf{[ILH.2']}, we conclude that $h=\<f,g\>$. Finally, by \textbf{[ILH.1']}, the evaluation map is a parametrized multilinear map. Furthermore, by construction, $(z\times\hat p)\ev=\T\ev\hat p$, which is equivalent to linearity in the last argument. Thus, $\ev$ is multilinear.
\end{proof}

Theorem~\ref{theorem:internal-linear-homs-classify-kleisli-multilinear-morphisms} establishes what class of morphisms is truly classified by an internal multilinear Hom. Before moving on, we shall express condition \textbf{[ILH.2']} only in terms of \textbf{[ILH.2]}, but in a different category. To this end, recall that the tangent bundle functor $\T$ of a tangent category comes with a monad structure, whose unit is the zero vector field $z_M\colon M\to\T M$ and whose multiplication is $\mu_M\colon\T\T M\to\T M$, given by $\<\T p_M,p_{\T M}\>s_M$~\cite[Proposition~3.4]{cockett:tangent-cats}. Furthermore, the canonical flip $c_M\colon\T\T M\to\T\T M$ turns $\T$ into a strong tangent monad, which consists of a monad equipped with a distributive law between the tangent bundle functor and the monad (see~\cite[Definition~19]{cockett:tangent-monads}).

In particular, we shall call a \textbf{colax tangent monad} a monad $S\colon\X\to\X$ equipped with a distributive law $\alpha_M\colon\T SM\to S\T M$, compatible with the monad and the tangent structure; furthermore, we shall assume that the endofunctor $S$ preserves the fundamental tangent limits\footnote{These are the $n$-fold pullback of the projection along itself together with the image of this $n$-fold pullback along $\T^n$, and the pullback diagram that establishes the universal property of the vertical lift.}.

\begin{lemma}
\label{lemma:kleisli-category-of-tangent-bundle-functor}
The Kleisli category of a colax tangent monad $(S,\alpha)$ becomes a tangent category whose tangent bundle functor $\T^\KL$ sends each object $M$ to $\T M$ and each $f\colon M\nto N$ to $\T^\KL f\colon\T M\nto\T N$, whose intepretation is given by:
\begin{align*}
&\SQ{\T^\KL f}\colon\T M\xrightarrow{\T\SQ{f}}\T SN\xrightarrow{\alpha_N}S\T N
\end{align*}
Furthermore, the inclusion functor $(\X,\TT)\to\KL[S,\alpha]$ is a strict tangent morphism that preserves the limits preserved by $S$.
\end{lemma}
\begin{proof}
It is straightforward to lift the structural natural transformations of the tangent structure $\TT$ on the base tangent category to the Kleisli category. For instance, the projection $p^\KL$ is defined by $\SQ{p^\KL}\=p\eta$. The inclusion functor $\X\to\KL[S]$ sends each object $M$ to itself, and each morphism $f\colon M\to N$ to $f^\KL$ defined by $\SQ{f^\KL}\=f\eta$, where $\eta$ is the unit of $S$, and it preserves the limits preserved by the endofunctor $S$ itself. From the functoriality of the inclusion functor, all the equational axioms hold for $\TT^\KL$ to be a tangent structure; from the preservation of the fundamental tangent limits, also the universal properties of the fundamental tangent limits hold. Finally, by construction, the inclusion functor becomes a strict tangent morphism that preserves the limits preserved by $S$.
\end{proof}

\begin{remark}
\label{remark:kleisli-tangent-category}
In~\cite[Proposition~3.5]{lanfranchi:lie-correspondence-tangent-cats}, the author already proved that the \emph{co}Kleisli category of a \emph{lax} tangent \emph{co}monad is a tangent category and that the inclusion functor is a strict tangent morphism. One could see Lemma~\ref{lemma:kleisli-category-of-tangent-bundle-functor} as the dual statement of this proposition. However, one should also keep in mind that the dual argument does not apply directly, since the opposite category of a tangent category fails in general to be a tangent category. In particular, for Lemma~\ref{lemma:kleisli-category-of-tangent-bundle-functor} to hold, we required the monad $S$ to preserve the fundamental tangent limits, a condition that is not required for lax tangent comonads.
\end{remark}

Since $(\T,c)$ is a strong tangent monad, $(\T,c)$ becomes a colax tangent monad; so, in particular, its Kleisli category is a tangent category. Furthermore, as proved in~\cite[Corollary~5.31]{lanfranchi:phd-thesis}, not only $(\T,c)$ is a strong tangent monad, but so is every $(\T^n,c^n)$, where $c^n$ is given by $c^n\=\T^nc\T^{n-1}c_\T\dots\T c_{\T^{n-1}}c_{\T^n}$.

\begin{proposition}
\label{proposition:kleisli-tangent-category-of-tangent-bundle-functor}
Given a (Cartesian) tangent category $(\X,\TT)$, for each $n$, the Kleisli category $\KL[\T^n,c^n]$ of $(\T^n,c^n)$ is a (Cartesian) tangent category and the inclusion functor $(\X,\TT\to\KL[\T^n,c^n]$ is a strict tangent morphism (that preserves Cartesian products). In particular, it preserves differential objects and multilinear maps.
\end{proposition}

We can now express condition \textbf{[ILH.2']} in a clearer form.

\begin{corollary}
\label{corollary:internal-linear-homs-classify-kleisli-multilinear-morphisms}
Given a list $\A_1\,\A_n,\B$ of differential objects, consider an object $\lnr[\A_1\,\A_n;\B]$ equipped with a map $\ev\colon\A_1\ntimes\A_n\times\lnr[\A_1\,\A_n;\B]\to\B$. Then $\lnr[\A_1\,\A_n;\B]$ is an internal multilinear Hom with evaluation map $\ev$ if and only if the following conditions are satisfied:
\begin{description}
\item[ILH.1'] $\ev$ is a $\lnr[\A_1\,\A_n;\B]$-parametrized multilinear map;

\item[ILH.2''] The morphism $\ev^\KL$ defined by $\SQ{\ev^\KL}=\ev z$, classifies multilinear maps in the Kleisli tangent category $\KL[\T,c]$, in that for every $M$-parametrized multilinear map $f\colon\A_1\ntimes\A_n\times M\nto\B$, there exists a unique morphism $\hat f\colon M\nto\lnr[\A_1\,\A_n;\B]$ that makes the following diagram
\begin{equation*}
\begin{tikzcd}
{\A_1\ntimes\A_n\times\lnr[\A_1\,\A_n;\B]} & \B \\
{\A_1\ntimes\A_n\times M}
\arrow["{\ev^\KL}"{inner sep=.8ex}, "\shortmid"{marking}, from=1-1, to=1-2]
\arrow["{\id\times\hat f}"{inner sep=.8ex}, "\shortmid"{marking}, from=2-1, to=1-1]
\arrow["f"'{inner sep=.8ex}, "\shortmid"{marking}, from=2-1, to=1-2]
\end{tikzcd}
\end{equation*}
commutative in $\KL[\T,c]$.
\end{description}
\end{corollary}
\begin{proof}
To prove this result, we shall simply rewrite condition \textbf{[ILH.2']} as \textbf{[ILH.2'']}. To this end, consider an $M$-parametrized multilinear map $f\colon\A_1\ntimes\A_n\times M\to\T\B$. This is clearly a morphism in the Kleisli category $\KL[\T,c]$. On the other hand, $\ev^\KL$ is defined by $\SQ{\ev^\KL}=\ev z$. Thus, by Kleisli composition, which we shall denote by $;_\KL$, for a given morphism $\hat f\colon M\to\T\lnr[\A_1\,\A_n;\B]$:
\begin{align*}
\SQ{(\id\times\hat f);_\KL\ev^\KL}&=~(\SQ{\id}\times\hat f)\T\SQ{\ev^\KL}\mu    \tag{Kleisli composition and products in the Kleisli}\\
&=~(z\times\hat f)\T\ev\T z\mu                                                  \tag{Identity in the Kleisli category and $\SQ{\ev^\KL}=\ev z$}\\
&=~(z\times\hat f)\T\ev                                                         \tag{Unitality: $\T z\mu=\id$}
\end{align*}
Therefore, $(\id\times\hat f);_\KL\ev^\KL=f$ in $\KL[\T,c]$ if and only if $(z\times\hat f)\T\ev=f$ in the base tangent category.
\end{proof}

We shall refer to the morphisms in the Kleisli category as \textbf{Kleisli morphisms}. We shall also refer to the $M$-parametrized multilinear morphisms in the Kleisli tangent category $\KL[\T^n,c^n]$ as \textbf{Kleisli $M$-parametrized multilinear morphisms} of order $n$.

\subsection{The algebraic structures of Linear Endomorphism Classifiers}
\label{subsection:ring-lie-algebra-lec}
For the story of this paper, instead of considering every internal multilinear Hom $\lnr[\A_1\,\A_n;\B]$, we shall focus on the special case $\lnr[\A,\A]$.

\begin{definition}
\label{definition:linear-endomorphism-classifier}
A differential object $\A$ admits the \textbf{Linear Endomorphism Classifier} provided that the internal multilinear Hom $\lnr_\A\=\lnr[\A,\A]$ exists. we will refer to $\lnr_\A$ as the \textbf{Linear Endomorphism Classifier} (\textbf{LEC}) of $\A$.
\end{definition}

In this section, we study the algebraic structures of LEC objects. We start with the most obvious one: every LEC is a rig object\footnote{Rigs are also known in the literature as semirings. Here, we shall not assume commutativity of the multiplication.}. Recall that a \textbf{differential rig object} consists of a differential object $\R$ equipped with a unit $\e\colon\*\to\R$ and a bilinear multiplication map $\m\colon\R\times\R\to\R$ which is associative and unital (see~\cite[Definition~5.24]{cockett:differential-equations}). We notice that such an object is equivalent to a monoid in the multicategory $\DO_\multi[\X,\TT]$.

\begin{proposition}
\label{proposition:rig-structure-on-lec}
The LEC object $\lnr_\A$ of a differential object $A$ is a differential rig object whose unit and bilinear multiplication $\m\colon\lnr_\A\times\lnr_\A\to\lnr_\A$ are respectively defined by:
\begin{equation*}
\begin{tikzcd}
{\A\times\lnr_\A} & \A \\
{\A\times\*}
\arrow["{\ev_\A}", from=1-1, to=1-2]
\arrow["{\id_\A\times\e}", from=2-1, to=1-1]
\arrow["{\pi_1}"', from=2-1, to=1-2]
\end{tikzcd}\qquad
\begin{tikzcd}
{\A\times\lnr_\A} & \A \\
{\A\times\lnr_\A\times\lnr_\A} & {\A\times\lnr_\A}
\arrow["{\ev_\A}", from=1-1, to=1-2]
\arrow["{\id_\A\times\m}", from=2-1, to=1-1]
\arrow["{\ev_\A\times\id_{\lnr_\A}}"', from=2-1, to=2-2]
\arrow["{\ev_\A}"', from=2-2, to=1-2]
\end{tikzcd}
\end{equation*}
\end{proposition}
\begin{proof}
By Proposition~\ref{proposition:exponential-objects}, $\lnr_\A=\lnr[\A,\A]$ is an exponential object in the multicategory $\DO_\multi[\X,\TT]$. Thus, it becomes a monoid object in the multicategory whose multiplication $\A^\A,\A^\A\to\A^\A$ is the uncurry of $(\ev_\A\times\A)\ev_\A\colon\A\times\A^\A\times\A^\A\to\A$ and whose unit $\text{\textvisiblespace}\to\A^\A$ is the uncurry of the identity $\id\colon\A\to\A$.
\end{proof}

One can turn any ring $R$ into a Lie algebra by considering the Lie bracket given by the commutator $[r,s]\=sr-rs$. Unfortunately, in order to define the commutator, one needs negatives, which are not assumed to exist in a rig. To solve this issue, we shall consider the maximal subring of a rig object.

\begin{definition}
\label{definition:maximal-differential-subring}
In a Cartesian tangent category, the \textbf{maximal differential subring} of a differential rig $\R$ consists of a differential object $\R^-$ together with a linear morphism $\varphi\colon\R^-\to\R\times\R$ that turns the following diagram
\begin{equation*}
\begin{tikzcd}
{\R^-} & {\R\times\R} & {\*} & \R
\arrow["\varphi", dashed, from=1-1, to=1-2]
\arrow[shift right=2, from=1-2, to=1-3]
\arrow["{+}", shift left=2, from=1-2, to=1-4]
\arrow["0"', shift right=2, from=1-3, to=1-4]
\end{tikzcd}
\end{equation*}
into an equalizer.
\end{definition}

The maximal differential ring of a differential rig becomes a differential ring that turns the inclusion morphism $\rho\colon\R^-\to\R$ into a linear rig homomorphism.

\begin{lemma}
\label{lemma:maximal-differential-subring}
The maximal differential ring $\R^-$ of a differential rig $\R$ comes with a unique differential ring structure that makes $\rho\=\varphi\pi_1\colon\R^-\to\R$ into a linear rig homomorphism.
\end{lemma}
\begin{proof}
This result applies more generally to rig objects and their maximal subring objects and is an easy exercise of categorical algebra. Concretely, to construct the negation morphism $-\colon\R^-\to\R^-$ one considers the morphism $\varphi\pi_2$, where $\varphi\colon\R^-\to\R\times\R$ denotes the equalizer morphism. Since $\varphi+=!0$, it follows that $\<!0,\varphi\pi_2\>+=\<!0\varphi\pi_1,\varphi\pi_2\>+=\<!0,\id\>\varphi+=!0$. Therefore, by the universal property of the equalizer, there exists a unique morphism $-\colon\R^-\to\R^-$ such that $-\varphi=\<!0,\varphi\pi_2\>$. One can show that $-$ defines the negation. We leave it to the reader to show that for each rig object $R$, the equalizer $R^-$ of $+$ and $0$ is in fact a ring object. Furthermore, it is also straightforward to see that the inclusion map $\rho$ becomes a rig homomorphism and that, by the universal property, the ring structure of $\R^-$ is uniquely determined by $\rho$ being a rig homomorphism. Finally, to conclude, we need to prove that $\R^-$ is a \emph{differential} ring. To this end, we shall use that $\R^-$ is a differential object and that $\varphi\colon\R^-\to\R\times\R$ is linear. Since $\varphi$ is linear, it follows that $\rho=\varphi\pi_1\colon\R^-\to\R$ is linear as well. Let $\m^-$ denote the multiplication of $\R^-$. Using that $\rho$ is a linear rig homomorphism, and that $\m$ is bilinear, it is not hard to prove that:
\begin{align*}
&(\id\times z)\T\m^-\hat p^-\rho=(\hat p^-\times\id)\m^-\rho &&(z\times\id)\T\m^-\hat p^-\rho=(\id\times\hat p^-)\m^-\rho
\end{align*}
Furthermore, using that $\varphi=\<\rho,-\rho\>$ since $-\varphi=\<!0,\varphi\pi_2\>$, it follows that $\rho$ is monic and therefore, $\m^-$ becomes bilinear.
\end{proof}

Every differential ring comes with a differential Lie algebra structure, induced by the commutator. For starters, let us recall the notion of differential Lie algebra.

\begin{definition}[{\cite[Definition~5.1]{lanfranchi:lie-correspondence-tangent-cats}}]
\label{definition:differential-lie-algebra}
In a Cartesian tangent category, a \textbf{differential Lie algebra} consists of a differential object $\g$ equipped with a morphism $[,]\colon\g\times\g\to\g$ called the \textbf{Lie bracket} of $\g$, satisfying the following conditions:
\begin{description}
\item[DLIE.1] $(\g,0,+,[,])$ is a Lie algebra, in that\footnote{For our convenience, we use pointwise notation to express antisymmetry and the Jacobian identity, as if $x,y$ and $z$ were internal variables of $\g$; however, it should be clear that these axioms can also be fully expressed diagrammatically.}:
\begin{align*}
&[x,y]+[y,x]=0              &&[[x,y],z]+[[z,x],y]+[[y,z],x]=0
\end{align*}

\item[DLIE.2] $[,]\colon\g\times\g\to\g$ is bilinear.
\end{description}
\end{definition}

\begin{proposition}
\label{proposition:LEC-lie-algebra}
The commutator $[,]$ of a differential ring $\R$ equips $\R$ with a differential Lie algebra structure. In particular, the maximal differential subring of a differential rig is a differential Lie algebra.
\end{proposition}
\begin{proof}
To prove that $[,]$ satisfies the antisymmetry identity, we shall use an internal language proof:
\begin{align*}
&[x,y]+[y,x]=\m^\op(x,y)-\m(x,y)+\m^\op(y,x)-\m(y,x)=\\
&\qquad=\m(y,x)-\m(x,y)+\m(x,y)-\m(y,x)=0
\end{align*}
Similarly, using the associativity of $\m$, one can prove the Jacobi identity with a similar technique. We leave the details to the reader. Finally, since $\m$ is bilinear, it is not hard to see that the linear combination $\m^\op-\m$ is also bilinear, proving that $\lnr_\A$ is in fact a differential Lie algebra.
\end{proof}

In particular, the maximal differential subring $\lnr_\A^-$ of the LEC object $\lnr_\A$ of a differential object $\A$ becomes a differential Lie algebra with Lie bracket given by the commutator.


\section{General linear groups in tangent categories}
\label{section:general-linear-groups}
In this section, we introduce general linear groups as objects that classify \emph{Kleisli parametrized linear automorphisms}, which are those morphisms $f\colon\A\times M\to\T^n\A$ which are linear in their first argument and that satisfy a suitable notion of invertibility. After introducing this concept, we prove that, as the name suggests, general linear groups are in fact group objects. We then establish that general linear groups can be constructed entirely from LEC objects $\lnr_\A$ via a tangent triple equalizer.

\subsection{General linear group objects}
\label{subsection:general-linear-groups}
Corollary~\ref{corollary:internal-linear-homs-classify-kleisli-multilinear-morphisms} establishes that the LEC object $\lnr_\A$ classifies Kleisli parametrized linear morphisms $f\colon\A\times M\to\T\A$. We shall define general linear group objects as those objects that classify those parametrized linear morphisms $f\colon\A\times M\to\T\A$ which are \emph{invertible}. The issue is that we do not expect any $f\colon\A\times M\to\T\A$ to be invertible since $\A$ and $\T\A$ will not be isomorphic objects. To find the correct notion of invertibility, we shall interpret a morphism $f\colon\A\times M\to\T\A$ as an \emph{endomorphism} in the appropriate category.

\begin{definition}
\label{definition:parametrized-linear-isomorphism}
Given two differential objects $\A$ and $\B$, an \textbf{$M$-parametrized linear isomorphism} consists of a morphism $f\colon\A\times M\to\B$ that satisfies the following conditions:
\begin{description}
\item[PLI.1] $f$ is an $M$-parametrized linear morphism, in that:
\begin{equation*}
\begin{tikzcd}
{\T(\A\times M)} & {\T\B} \\
{\A\times M} & \B
\arrow["{\T f}", from=1-1, to=1-2]
\arrow["{\lambda_\A\times z_M}", from=2-1, to=1-1]
\arrow["f"', from=2-1, to=2-2]
\arrow["{\lambda_\B}"', from=2-2, to=1-2]
\end{tikzcd}
\end{equation*}

\item[PLI.2] There exists a morphism $f^{-M}\colon\B\times M\to\A$, called \textbf{$M$-parametrized inverse}, that makes the following diagrams:
\begin{equation*}
\begin{tikzcd}
& {\A\times M\times M} \\
{\A\times M} & {\B\times M} \\
& \A
\arrow["{f\times\id_M}", from=1-2, to=2-2]
\arrow["{\id_\A\times\Delta_M}", from=2-1, to=1-2]
\arrow["{\pi_1}"', from=2-1, to=3-2]
\arrow["{f^{-M}}", from=2-2, to=3-2]
\end{tikzcd}\qquad
\begin{tikzcd}
{\B\times M\times M} & \\
{\A\times M} & {\B\times M} \\
\B
\arrow["{f^{-M}\times\id_M}"', from=1-1, to=2-1]
\arrow["f"', from=2-1, to=3-1]
\arrow["{\id_\B\times\Delta_M}"', from=2-2, to=1-1]
\arrow["{\pi_1}", from=2-2, to=3-1]
\end{tikzcd}
\end{equation*}
commutative.
\end{description}
When $\A=\B$, we will refer to such an $f$ as a \textbf{parametrized linear automorphism} of $\A$.
\end{definition}

In Lemma~\ref{lemma:multilinearity-fibrationally}, we established that $M$-parametrized linear morphisms $f\colon\A\times M\to\B$ can be equivalently characterized as \emph{linear} morphisms in $\Simple_M[\X,\TT]$. The next lemma establishes a similar characterization for parametrized linear isomorphisms.

\begin{lemma}
\label{lemma:parametrized-linear-isomorphisms}
A morphism $f\colon\A\times M\to\B$ is a parametrized linear isomorphism if and only if it lifts to a linear isomorphism $f\colon\A\nto\B$ in $\Simple_M[\X,\TT]$. Furthermore, the $M$-parametrized inverse $f^{-M}$ of $f$ is unique.
\end{lemma}
\begin{proof}
From Lemma~\ref{lemma:multilinearity-fibrationally}, $f\colon\A\times M\to\B$ is an $M$-parametrized linear map if and only if $f\colon\A\nto\B$ is linear in $\Simple_M[\X,\TT]$. It is also easy to see that the conditions expressed by the diagrams of \textbf{[PLI.2]} establish that $f^{-M}$ is an inverse of $f$ in $\Simple_M[\X,\TT]$. Finally, since inverses of morphisms are always unique, it follows that parametrized inverses are also unique.
\end{proof}

For our goal, we need to impose an invertibility condition on morphisms of type $f\colon\A\times M\to\T\A$. However, if $f$ is an $M$-parametrized linear isomorphism, then $\A\cong\T\A$, which is not the case. Instead, we shall reinterpret this class of morphisms as endomorphisms in the Kleisli tangent category of $(\T,c)$.

\begin{definition}
\label{definition:kleisli-PLAs}
An \textbf{Kleisli $M$-parametrized linear isomorphism} of order $n$ consists of a morphism $f\colon\A\times M\to\T^n\B$ subject to the following conditions:
\begin{description}
\item[KPLI.1] $f$ is an $M$-parametrized linear morphism;

\item[KPLI.2] There exists a morphism $f^{-M}\colon\B\times M\to\T^n\A$, called \textbf{Kleisli $M$-parametrized inverse} of $f$, that satisfies the following conditions:
\begin{equation*}
\begin{tikzcd}
{\A\times M} & {\A\times M\times M} & {\T^n(\A\times M)} \\
&& {\T^n\T^n\A} \\
\A && {\T^n\A}
\arrow["{\id\times\Delta}", from=1-1, to=1-2]
\arrow["{\pi_1}"', from=1-1, to=3-1]
\arrow["{f\times z^n}", from=1-2, to=1-3]
\arrow["{\T^nf^{-M}}", from=1-3, to=2-3]
\arrow["{\mu^n}", from=2-3, to=3-3]
\arrow["{z^n}"', from=3-1, to=3-3]
\end{tikzcd}
\end{equation*}
\begin{equation*}
\begin{tikzcd}
{\A\times M} & {\A\times M\times M} & {\T^n(\A\times M)} \\
&& {\T^n\T^n\A} \\
\A && {\T^n\A}
\arrow["{\id\times\Delta}", from=1-1, to=1-2]
\arrow["{\pi_1}"', from=1-1, to=3-1]
\arrow["{f^{-M}\times z^n}", from=1-2, to=1-3]
\arrow["{\T^nf}", from=1-3, to=2-3]
\arrow["{\mu^n}", from=2-3, to=3-3]
\arrow["{z^n}"', from=3-1, to=3-3]
\end{tikzcd}
\end{equation*}
\end{description}
where $z^n$ and $\mu^n$ denote the unit and the multiplication of the colax tangent monad $(\T^n,c^n)$. When $\A=\B$, we will refer to such a map as a \textbf{Kleisli Parametrized Linear Automorphism}.
\end{definition}

\begin{lemma}
\label{lemma:kleisli-PLI}
An $M$-parametrized linear morphism $f\colon\A\times M\to\T^n\B$ is a Kleisli $M$-parametrized linear isomorphism if and only if $f$ lifts to a linear isomorphism $f\colon\A\nto\B$ in the Kleisli tangent category $\KL[\T_M^n,c^n_M]$ of the colax tangent monad $(\T_M^n,c_M^n)$, where $\T^n_M$ denotes the iterate of the tangent bundle functor $\T_M$ of $\Simple_M[\X,\TT]$. Moreover, the Kleisli parametrized inverse $f^{-M}$ of $f$ is unique.
\end{lemma}
\begin{proof}
That $f$ is a linear isomorphism in the Kleisli tangent category follows directly from the definition. Furthermore, uniqueness of the Kleisli parametrized inverse follows from the uniqueness of the inverse of an isomorphism.
\end{proof}

We can now define general linear group objects.

\begin{definition}
\label{definition:general-linear-group}
The \textbf{general linear group object} of $\A$ consists of an object $\GL_\A$ equipped with a map $\Ev\colon\A\times\GL_\A\to\A$, called the \textbf{evaluation map}, that satisfies the following conditions:
\begin{description}
\item[GL.1] $\Ev$ is a $\GL_\A$-parametrized linear automorphism.

\item[GL.2] $\Ev$ classifies Kleisli parametrized linear automorphisms of $\A$, in that if $f\colon\A\times M\to\T^n\A$ is a Kleisli $M$-parametrized linear automorphism of $\A$, there exists a unique $f\colon M\to\T^n\GL_\A$ that makes the following diagram
\begin{equation*}
\begin{tikzcd}
{\T^n(\A\times\GL_\A)} & {\T^n\A} \\
{\A\times M}
\arrow["{\T^n\Ev}", from=1-1, to=1-2]
\arrow["{z^n\times\hat f}", from=2-1, to=1-1]
\arrow["f"', from=2-1, to=1-2]
\end{tikzcd}
\end{equation*}
commutative.
\end{description}
\end{definition}

\begin{example}
\label{example:GL-trivial}
Every category is trivially a tangent category whose tangent bundle functor and structural natural transformations are just the identity. When the category is Cartesian, it also becomes a Cartesian tangent category. However, the only differential object is the terminal object $\*$. It is also immediate to see that $\*$ always admits the general linear group, $\GL_\*$ being $\*$ and the evaluation map $\Ev_\*\colon\*\times\*\to\*$ being the unique terminal map.
\end{example}

\begin{example}
\label{example:GL-sman}
In the tangent category $\SMan$ of finitely-dimensional smooth manifolds, differential objects are precisely Euclidean spaces $\Real^n$ (\cite[Example~3.3(i)]{cockett:differential-bundles}). Every differential object admits the general linear group, whose $\GL$-object is precisely $\GL_n(\Real)$ and the evaluation map $\Real^n\times\GL_n(\Real)\to\Real^n$ sends each vector $v\in\Real^n$ and each linear automorphism $l\colon\Real^n\to\Real^n$ to $l(v)$.
\end{example}

\begin{example}
\label{example:GL-affine}
Affine schemes over a commutative and unital ring $R$ form a Cartesian tangent category $\Aff_R$ whose tangent bundle functor sends an affine scheme with coordinate ring $A$ to $\T^\Omega A\=\Sym_A\Omega_A$, that is, the symmetric algebra of the module of K\"ahler differentials of $A$. See~\cite{cruttwell:algebraic-geometry} for details. By~\cite[Corollary~4.18 and Remark~4.20]{cruttwell:algebraic-geometry}, there is an equivalence of categories:
\begin{align*}
&\DO[\Aff_R]\simeq\Mod_R^\op
\end{align*}
In particular, every differential object $\A$ of $\Aff_R$ is the free symmetric $R$-algebra $\A\cong\Sym_RM$ of an $R$-module $M$. For simplicity, we consider the case when $R$ is the field $\Real$ of real numbers. Consider a differential object $\A$ whose $\Real$-module $M$ is finitely generated by $n$ generators. Since $\Real$ is a field, $M$ becomes a vector space of dimension $n$.

Now, define $\GL_\A(-)$ as the functor that sends an $\Real$-algebra $A\in\Aff_\Real$ to the group $\GL_\A(A)$ of invertible $(n\times n)$-matrices with coefficients in $A$. This functor is represented by the $\Real$-algebra defined as follows:
\begin{align*}
&\GL_\A\=\Real[\{x_{i,j}\}_{i,j=1}^n][\det(x_{i,j})]^{-1}
\end{align*}
that is, the localization of the $\Real$-algebra freely-generated by the $n^2$-variables $x_{i,j}$, localized to invert the determinant polynomial $\det(x_{i,j})$ of the matrix $(x_{i,j})$. To see this, consider a morphism of $\Aff_\Real$ from $A\to\GL_\A$. This corresponds to a morphism of $\Real$-algebras $\GL_\A\to A$ which is fully determined by its values $a_{i,j}$ at each variable $x_{i,j}$. However, since it is an algebra homomorphism, $\det(a_{i,j})$ must be invertible in $A$. That is, the matrix $(a_{i,j})$ is invertible, thus, $(a_{i,j})\in\GL_\A(A)$.

The object $\GL_\A$ turns out to be the general linear group of $\A$.
\end{example}

To support our definition of general linear groups, we begin by showing that $\GL_\A$ classifies linear automorphisms of $\A$.

\begin{lemma}
\label{lemma:gl-classifies-linear-automorphisms}
The general linear group object $\GL_\A$ of a differential object $\A$ satisfies the following universal property:
\begin{description}
\item[GL.2'] For every $M$-parametrized linear automorphism $f\colon\A\times M\to\A$, there exists exactly one morphism $\hat f\colon M\to\GL_\A$ that makes the following diagram
\begin{equation*}
\begin{tikzcd}
{\A\times\GL_\A} & \A \\
{\A\times M}
\arrow["\Ev", from=1-1, to=1-2]
\arrow["{\id_\A\times\hat f}", from=2-1, to=1-1]
\arrow["f"', from=2-1, to=1-2]
\end{tikzcd}
\end{equation*}
commutative.
\end{description}
\end{lemma}
\begin{proof}
Consider an $M$-parametrized linear automorphism $f\colon\A\times M\to\A$. We shall prove that $fz$ becomes a Kleisli $M$-parametrized linear automorphism of $\A$. For starters, since $z\T\lambda c=\lambda z_\T c=\lambda\T z$, it follows that $f$ is linear in the first argument. To prove that $fz$ is a Kleisli parametrized isomorphism, define $g\colon\A\times M\to\T\A$ as $f^{-M}z$, where $f^{-M}$ denotes the parametrized inverse of $f$. We shall compute:
\begin{align*}
(\id\times\Delta)(fz\times\id)(\id\times z)\T f^{-M}\T z\mu&=~(\id\times\Delta)(fz\times\id)(\id\times z)\T f^{-M}                      \tag{Unitality of the monad}\\
&=~(\id\times\Delta)(f\times\id)f^{-M}z             \tag{Naturality of $z$}\\
&=~(\id\times\Delta)\pi_1z                          \tag{$f^{-M}$ parametrized inverse of $f$}\\
&=~\pi_1z
\end{align*}
Similarly, we can also prove that $(\id\times\Delta)(f^{-M}z\times\id)(\id\times z)\T f\T z\mu=\pi_1z$. Therefore, $fz$ is a Kleisli $M$-parametrized linear automorphism of $\A$. Therefore, by the universal property \textbf{[GL.2]}, there exists a unique morphism $h\colon M\to\T\GL_\A$. Let $\hat f\colon M\to\GL_\A$, simply be $\hat f\=hp$. We shall prove that $\hat f$ satisfies $(\id\times\hat f)\Ev=f$:
\begin{align*}
(\id\times\hat f)\Ev&=~(zp\times hp)\Ev     \tag{Definition of $\hat f$ and $zp=\id$}\\
&=~(z\times h)\T\Ev p                       \tag{Naturality of $p$}\\
&=~fzp                                      \tag{Definition of $h$}\\
&=~f                                        \tag{$zp=\id$}
\end{align*}
Now, let us prove that $\hat f$ is unique. Let $t$ satisfy $(\id\times t)\Ev=f$. Thus:
\begin{align*}
&(z\times h)\T\Ev=fz=(\id\times t)\Ev z=(z\times tz)\T\Ev
\end{align*}
Thus, $tz=h$. Thus, $t=tzp=hp=\hat f$. This proves that $\hat f$ is unique.
\end{proof}

The universal property \textbf{[GL.2]} implies the universal property \textbf{[GL.2']} of Lemma~\ref{lemma:gl-classifies-linear-automorphisms}; however, it is strictly stronger, since the class of Kleisli parametrized automorphisms is bigger than the one of parametrized automorphisms. Although \textbf{[GL.2]} is strictly stronger than \textbf{[GL.2']}, we shall express \textbf{[GL.2]} entirely in terms of \textbf{[GL.2']}, using that Kleisli parametrized automorphisms are in fact Kleisli morphisms.

\begin{lemma}
\label{lemma:general-linear-groups-kleisli}
Consider a differential object $\A$ and an object $\GL_\A$ together with a map $\Ev\colon\A\times\GL_\A\to\A$. Thus, $\GL_\A$ is a general linear group for $\A$ with evaluation map $\Ev$ if and only if it satisfies the following conditions:
\begin{description}
\item[GL.1] $\Ev$ is a $\GL_\A$-parametrized linear automorphism of $\A$;

\item[GL.2''] For every $n>0$ and every parametrized linear automorphism $f\colon\A\times M\nto\A$ in $\KL[\T^n,c^n]$, there exists a unique morphism $\hat f\colon M\nto\GL_\A$ of $\KL[\T^n,c^n]$ that makes the following diagram
\begin{equation*}
\begin{tikzcd}
{\A\times\GL_\A} & \A \\
{\A\times M}
\arrow["{\Ev^\KL}"{inner sep=.8ex}, "\shortmid"{marking}, from=1-1, to=1-2]
\arrow["{\id_\A\times\hat f}"{inner sep=.8ex}, "\shortmid"{marking}, from=2-1, to=1-1]
\arrow["f"'{inner sep=.8ex}, "\shortmid"{marking}, from=2-1, to=1-2]
\end{tikzcd}
\end{equation*}
commutative in $\KL[\T^n,c^n]$, where $\SQ{\Ev^\KL}\=\Ev z^n$.
\end{description}
\end{lemma}
\begin{proof}
The argument for this statement is fairly similar to the one of~\ref{corollary:internal-linear-homs-classify-kleisli-multilinear-morphisms} as it consists entirely of rewriting condition \textbf{[GL.2]} as a condition in the Kleisli category of $\T^n$.
\end{proof}

By Lemma~\ref{lemma:gl-classifies-linear-automorphisms}, we can prove that $\GL_\A$ is a group object.

\begin{proposition}
\label{proposition:gl-group-object}
The general linear group object $\GL_\A$ of $\A$ is a group object whose unit, multiplication, and inverse maps are uniquely defined by the following commutative diagrams:
\begin{equation*}
\begin{tikzcd}
{\A\times\GL_\A} & \A \\
{\A\times\*}
\arrow["\Ev", from=1-1, to=1-2]
\arrow["{\id_\A\times\e}", from=2-1, to=1-1]
\arrow["{\pi_1}"', from=2-1, to=1-2]
\end{tikzcd}\qquad
\begin{tikzcd}
{\A\times\GL_\A} & \A \\
{\A\times\GL_\A\times\GL_\A} & {\A\times\GL_\A}
\arrow["\Ev", from=1-1, to=1-2]
\arrow["{\id_\A\times\m}", from=2-1, to=1-1]
\arrow["{\Ev\times\id}"', from=2-1, to=2-2]
\arrow["\Ev"', from=2-2, to=1-2]
\end{tikzcd}\qquad
\begin{tikzcd}
{\A\times\GL_\A} & \A \\
{\A\times\GL_\A}
\arrow["\Ev", from=1-1, to=1-2]
\arrow["{\id\times\i}", from=2-1, to=1-1]
\arrow["{\Ev^{-\GL_\A}}"', from=2-1, to=1-2]
\end{tikzcd}
\end{equation*}
where $\Ev^{-\GL_\A}$ denotes the parametrized inverse of $\Ev$.
\end{proposition}
\begin{proof}
One can directly check the equations that define a group object; however, we decide to give a more structural perspective on this result. We begin by constructing a functor $\GL_\A(-)\colon\X\to\Gr^\op$, where $\Gr$ is the category of groups. Let $\GL_\A(M)$ be the set of all $M$-parametrized linear automorphisms $f\colon\A\times M\to\A$. By Lemma~\ref{lemma:parametrized-linear-isomorphisms}, PLAs are linear isomorphisms in $\Simple_M[\X,\TT]$, thus, $\GL_\A(M)$ becomes a group object. Moreover, by Proposition~\ref{proposition:simple-tangent-fibration}, since $\Simple[\X,\TT]\to(\X,\TT)$ is a tangent fibration, the substitution functors $g^\*\colon\Simple_N[\X,\TT]\to\Simple_M[\X,\TT]$, induced by $g\colon M\to N$ are strong tangent morphisms. So, they preserve linear isomorphisms and, by functoriality, they restrict to group homomorphisms. Thus, we can define $\GL_\A(g)\=g^\*\colon\GL_\A(N)\to\GL_\A(M)$. It turns out that the universal property \textbf{[GL.2']}, established by Lemma~\ref{lemma:gl-classifies-linear-automorphisms} is equivalent to the representability of the functor $\GL_\A(-)$. We leave it to the reader to prove this statement. By the Yoneda lemma, the group structure of each $\GL_\A(M)$ induces a group structure on the representing object, which is precisely $\GL_\A$.
\end{proof}

\subsection{General linear groups as tangent equalizers}
\label{subsection:equalizers}
We now consider the case when $\A$ is a differential object which admits the LEC object $\lnr_\A$. We shall prove that one can construct the general linear group $\GL_\A$ of $\A$ out of $\lnr_\A$ via a \emph{tangent} triple equalizer. We begin with a technical lemma. Recall that a \textbf{tangent limit}, as introduced in~\cite[Definition~1.1.3]{macadam:vector-bundles}, consists of a limit cone $\gamma_c\colon L\to D(c)$ of a diagram $D\colon I\to\X$ which is preserved by all iterates of the tangent bundle functor, in that $\T^n\gamma_c\colon\T^nL\to\T D(c)$ is a limit cone for the diagram $\T^n\o D\colon I\to\X$. We shall also call a tangent limit of a particular shape, ``tangent [the name of that limit shape]'', e.g., a \textbf{tangent pullback} is a pullback diagram whose underlying cone is a tangent limit cone.

\begin{lemma}
\label{lemma:tangent-limits-kleisli}
Consider a diagram $D\colon I\to\X$ and a tangent structure $\TT$ on $\X$. A limit cone $\gamma_c\colon L\to D(c)$ for $D$ is a tangent limit if and only if, for every $n>0$, the inclusion functor $(\X,\TT)\to\KL[\T^n,c^n]$ preserves the limit cone $\gamma_c$.
\end{lemma}
\begin{proof}
This follows immediately from the fact that a limit cone $\gamma_c\colon L\to D(c)$ is preserved by a monad $S$ if and only if the inclusion functor $\X\to\KL[S]$ preserves such a limit cone.
\end{proof}

\begin{theorem}
\label{theorem:triple-equalizer}
Consider a differential object $\A$ and assume that the LEC object $\lnr_\A$ of $\A$ exists. Then, $\A$ admits the general linear group $\GL_\A$ if and only if the triple equalizer
\begin{equation*}
\begin{tikzcd}
{\GL_\A} & {\lnr_\A\times\lnr_\A} & {\*} & {\lnr_\A}
\arrow[dashed, from=1-1, to=1-2]
\arrow[from=1-2, to=1-3]
\arrow["\m", curve={height=-24pt}, from=1-2, to=1-4]
\arrow["{\m^\op}"', curve={height=24pt}, from=1-2, to=1-4]
\arrow["\e", from=1-3, to=1-4]
\end{tikzcd}
\end{equation*}
exists and is a tangent triple equalizer.
\end{theorem}
\begin{proof}
We begin by showing that the triple equalizer exists if and only if there exists an object $\GL_\A$ together with a map $\Ev\colon\A\times\GL_\A\to\A$ that satisfies \textbf{[GL.1]} and \textbf{[GL.2']} of Lemma~\ref{lemma:gl-classifies-linear-automorphisms}. Using Lemma~\ref{lemma:general-linear-groups-kleisli}, we later prove that such a triple equalizer is a tangent triple equalizer if and only if $\GL_\A$ is a general linear group object with evaluation map $\Ev$.

For starters, let us assume that $\GL_\A$ and $\Ev$ satisfy \textbf{[GL.1]} and \textbf{[GL.2']} and let us prove that $\GL_\A$ is in fact a triple equalizer of $\m$, $!\e$, and $\m^\op$. To this end, we note that, since by \textbf{[GL.1]} $\Ev$ is a parametrized linear morphism, there exists a unique morphism $\theta\colon\GL_\A\to\lnr_\A$ that makes the following diagram
\begin{equation*}
\begin{tikzcd}
{\A\times\lnr_\A} & \A \\
{\A\times\GL_\A}
\arrow["\ev", from=1-1, to=1-2]
\arrow["{\id_\A\times\theta}", from=2-1, to=1-1]
\arrow["\Ev"', from=2-1, to=1-2]
\end{tikzcd}
\end{equation*}
commutative. We shall leave it to the reader to prove that $\theta$ is in fact a monoid homomorphism. Thus, we shall define $\varphi\colon\GL_\A\to\lnr_\A\times\lnr_\A$ as $\varphi\=\<\theta,\i\theta\>$, where $\i\colon\GL_\A\to\GL_\A$ denotes the inverse map of the group object $\GL_\A$\footnote{Note that in Proposition~\ref{proposition:gl-group-object} we only used \textbf{[GL.2']} to prove that $\GL_\A$ is a group object.}. We shall prove that $\varphi$ equalizes $\m$, $!\e$, and $\m^\op$:
\begin{align*}
\varphi\m &=~\<\theta,\i\theta\>\m\\
&=~\<\id,\i\>\m\theta               \tag{$\theta$ preserves the multiplication}\\
&=~!\e\theta                        \tag{Inverse law of $\GL_\A$}\\
&=~!\e                              \tag{$\theta$ preserves the unit}
\end{align*}
Similarly, one can also prove that $\varphi\m^\op=!\e$. Now, consider a morphism $f=\<f_1,f_2\>\colon M\to\lnr_\A\times\lnr_\A$ that equalizes $\m$, $!\e$, and $\m^\op$. Consider also the two morphisms:
\begin{align*}
&\tilde f_k\colon\A\times M\xrightarrow{\id\times f_k}\A\times\lnr_\A\xrightarrow{\ev}\A,    &&\text{for }k=1,2.
\end{align*}
It is easy to see that $\tilde f_1$ and $\tilde f_2$ are parametrized linear endomorphisms of $\A$; we shall also prove that they are parametrized inverses to each other:
\begin{align*}
(\id\times\Delta)(\tilde f_1\times\id)\tilde f_2&=~(\id\times\<f_1,f_2\>)(\ev\times\id)\ev   \tag{Definitions of $\tilde f_1$ and $\tilde f_2$}\\
&=~(\id\times f)(\id\times\m)\ev                                    \tag{Definition of $f$ and of $\m$}\\
&=~(\id\times!\e)\ev                                                \tag{$f\m=!\e$}\\
&=~\pi_1                                                            \tag{Definition of $\e$}
\end{align*}
Similarly, using that $f\m^\op=!\e$, one can also show that $(\id\times\Delta)(\tilde f_2\times\id)\tilde f_1=\pi_1$, proving that $\tilde f_1$ is a parametrized linear automorphism of $\A$. Therefore, by \textbf{[GL.2']}, there exists a unique morphism $\hat f\colon M\to\GL_\A$ such that $(\id\times\hat f)\Ev=\tilde f_1$. However, we can also compute that:
\begin{align*}
&(\id\times f_1)\ev=\tilde f_1=(\id\times\hat f)\Ev=(\id\times(\hat f\theta)\ev
\end{align*}
Thus, by \textbf{[GL.2']}, $\hat f\theta=f_1$. Moreover, there exists also a unique morphism $\hat f^{-M}\colon M\to\GL_\A$ such that $(\id\times\hat f^{-M})\Ev=\tilde f_2$, since also $\tilde f_2$ becomes a parametrized linear automorphism of $\A$. However:
\begin{align*}
&(\id\times f_2)\ev=\tilde f_2=(\id\times\hat f^{-M})\Ev= (\id\times(\hat f^{-M}\theta))\ev
\end{align*}
Thus, $f_2=\hat f^{-M}\theta$. However, using the inverse law, $(\id\times(\hat f\i)\theta)\ev$ is a parametrized inverse for $\tilde f_1$. Thus, by uniqueness of the parametrized inverse, $(\id\times(\hat f\i)\theta)\ev=\tilde f_2$; thus, $f_2=\hat f^{-M}\theta=(\hat f\i)\theta$. Therefore, $\hat f\varphi=\hat f\<\theta,\i\theta\>=\<f_1,f_2\>=f$. Now, let us consider another morphism $g\colon M\to\GL_\A$ that satisfies $g\varphi=f$. Thus, $g\theta=f_1=\hat f\theta$. Thus, $(\id\times g)\Ev=(\id\times g\theta)\ev=(\id\times\hat f\theta)\ev=(\id\times\hat f)\Ev$. Thus, by \textbf{[GL.2']}, $g=\hat f$. This proves that $\GL_\A$ is the triple equalizer of $\m$, $!\e$, and $\m^\op$.

We shall now prove the converse; so, we shall assume that $\GL_\A$ is an object that comes with a map $\varphi\=\<\varphi_1,\varphi_2\>\colon\GL_\A\to\lnr_\A\times\lnr_\A$ that makes it into the triple equalizer of $\m$, $!\e$, and $\m^\op$. We shall prove that, under this assumption, $\GL_\A$ satisfies \textbf{[GL.1]} and \textbf{[GL.2']}. We begin by constructing the evaluation map $\Ev$ as follows:
\begin{align*}
&\Ev\colon\A\times\GL_\A\xrightarrow{\id\times\varphi_1}\A\times\lnr_\A\xrightarrow{\ev}\A
\end{align*}
Let also:
\begin{align*}
&\Ev^{-\GL_\A}\colon\A\times\GL_\A\xrightarrow{\id\times\varphi_2}\A\times\lnr_\A\xrightarrow{\ev}\A
\end{align*}
We already proved before that, if $f=\<f_1,f_2\>$ equalizes $\m$, $!\e$, and $\m^\op$, then $(\id\times f_1)\ev$ and $(\id\times f_2)\ev$ are in fact parametrized linear automorphisms and parametrized inverses of each other. Therefore, by construction, so are $\Ev$ and $\Ev^{-\GL_\A}$. Thus, $\Ev$ satisfies \textbf{[GL.1]}. To prove \textbf{[GL.2']}, consider a parametrized linear automorphism $f\colon\A\times M\to\A$ whose parametrized inverse is denoted by $f^{-M}$. Using the universal property \textbf{[ILH.2]}, we find two morphisms $\hat f,\hat f^{-M}\colon M\to\lnr_\A$ satisfying $(\id\times\hat f)\ev=f$ and $(\id\times\hat f^{-M})\ev$. We shall prove that $\<\hat f,\hat f^{-M}\>$ equalizes $\m$, $!\e$, and $\m^\op$:
\begin{align*}
(\id\times(\<\hat f,\hat f^{-M}\>\m)\ev&=~(\id\times\<\hat f,\hat f^{-M}\>)(\ev\times\id)\ev    \tag{Definition of $\m$}\\
&=~(\id\times\Delta)((\id\times\hat f)\ev\times\id)(\id\times\hat f^{-M})\ev\\
&=~(\id\times\Delta)(f\times\id)f^{-M}              \tag{Definitions of $\hat f$ and $\hat f^{-M}$}\\
&=~\pi_1                                            \tag{$f^{-M}$ parametrized inverse of $f$}\\
&=~(\id\times!\e)\ev                                \tag{Definition of $\e$}
\end{align*}
Thus, by \textbf{[ILH.2]}, $\<\hat f,\hat f^{-M}\>\m=!\e$. Similarly, one can also show that $\<\hat f,\hat f^{-M}\>\m^\op=!\e$. Thus, by the universal property of the triple equalizer, there exists a unique morphism $\tilde f\colon M\to\GL_\A$ such that $\tilde f\varphi=\<\hat f,\hat f^{-M}\>$. Thus, $\tilde\varphi_1=\hat f$. Therefore, $(\id\times\tilde f)\Ev=(\id\times\tilde f\varphi_1)\ev=(\id\times\hat f)\ev=f$. We leave it to the reader to prove that, in fact, $\tilde f$ is the unique morphism that satisfies this equation.

So far, we proved that $\GL_\A$ is the triple equalizer of $\m$, $!\e$, and $\m^\op$ if and only if it comes with a map $\Ev$ that satisfies \textbf{[GL.1]} and \textbf{[GL.2']}. To prove the theorem, we shall make use of Lemma~\ref{lemma:general-linear-groups-kleisli}. According to this lemma, $\GL_\A$ is a general linear group with evaluation map given by $\Ev$ if and only if $\Ev$ satisfies \textbf{[GL.1]} and it satisfies \textbf{[GL.2']} in each Kleisli category $\KL[\T^n,c^n]$. Therefore, $\GL_\A$ is a general linear group if and only if, in each Kleisli category $\KL[\T^n,c^n]$, $\GL_\A$ is the triple equalizer of $\m$, $!\e$, and $\m^\op$. However, by Lemma~\ref{lemma:tangent-limits-kleisli}, a limit cone is a tangent limit if and only if it is a limit cone in each Kleisli category $\KL[\T^n,c^n]$. Therefore, $\GL_\A$ is a general linear group if and only if $\GL_\A$ is the \emph{tangent} triple equalizer of $\m$, $!\e$, and $\m^\op$.
\end{proof}

Conceptually, Theorem~\ref{theorem:triple-equalizer} establishes that the general linear group $\GL_\A$ is the maximal subgroup of $\lnr_\A$ and that the embedding is compatible with the tangent bundle functor.


\section{General linear Lie algebras}
\label{section:general-linear-algebras}
In this section, we study the Lie algebra of general linear groups. We shall prove that, in fact, such a Lie algebra $\gl_\A$ coincides with $\lnr_\A^-$, which is the differential Lie algebra of the maximal subring of the LEC object $\lnr_\A$, whose Lie bracket is given by the commutator. In particular, when $\A$ has negatives, this makes $\gl_\A$ isomorphic to $\lnr_\A$ itself, recapturing the classical result of differential geometry.

\subsection{The Lie algebra of a general linear group object}
\label{subsection:general-linear-algebra}
Recall that, as explained in the introduction, a \textbf{Lie group} in a Cartesian tangent category consists of a group object $G$ which admits the tangent space at the unit, in that the tangent pullback $\T_\e G$ of the tangent bundle $p_G\colon\T G\to G$ along the unit map $\e\colon\*\to G$ exists~\cite[Definition~2.13]{lanfranchi:lie-correspondence-tangent-cats}.

Recall also that the Lie functor $\g(-)$ of $G$ is the functor that sends each $M$ to the Lie algebra $\g(M)$ of $M$-parametrized left-invariant vector fields ($M$-\PLIVFs for short) of $G$. The Lie algebra of a Lie group is the representing object $\g$ of the Lie functor $\g(-)$ of $G$, which coincides with the tangent space $\T_\e G$. By~\cite[Proposition~4.3]{lanfranchi:lie-correspondence-tangent-cats}, this mounts to say that there is a \PLIVF $\omega\colon G\times\g\to\T G$ which is universal, in that for any $M$-\PLIVF $X\colon G\times M\to\T G$ of $G$ there exists a unique morphism $\hat X\colon M\to\g$ such that $X=(\id\times\hat X)\omega$. The universal \PLIVF of a Lie group is determined by the inclusion $\iota_\e\colon\T_\e G\to\T G$.

\begin{lemma}[{\cite[Lemma~4.5]{lanfranchi:lie-correspondence-tangent-cats}}]
\label{lemma:omega-vs-iota}
For a Lie group $G$, the following identity holds:
\begin{align*}
&\omega_G=(z\times\iota_\e)\T\m_G
\end{align*}
where $\omega_G\colon G\times\T_\e G\to\T G$ denotes the universal \PLIVF of $G$ and $\iota_\e\colon\T_\e G\to\T G$ is the inclusion of the tangent space in the tangent bundle of $G$.
\end{lemma}

Our goal is to compare the differential Lie algebra $\gl_\A$ of the general linear group $\GL_\A$ with the differential Lie algebra $\lnr_\A^-$ of Proposition~\ref{proposition:LEC-lie-algebra}. To make this comparison, we shall first re-interpret the existence of the Lie algebra object $\lnr_\A^-$ via the representability of the functor
\begin{align*}
&\lnr_\A^-(-)\colon\X\to\Set^\op
\end{align*}
that sends each object $M$ to the set $\lnr_\A^-(M)$ of $M$-parametrized linear endomorphisms $f\colon\A\times M\to\A$ of $\A$ for which there exists a morphism $-f$ such that, their sum $f+(-f)\=\<f,-f\>+$ vanishes, that is, $f+(-f)=!0$.

Another approach to construct the functor $\lnr_\A^-(-)$ consists of considering first the functor
\begin{align*}
&\lnr_\A(-)\colon\X\to\Set^\op  &&M\mapsto\lnr_\A(M)\=\DO_\multi(\Simple_M[\X,\TT])(\A,\A)
\end{align*}
that sends each $M$ to the set of $M$-parametrized linear endomorphisms. When the LEC object of $\A$ exists, this functor is representable and represented by the LEC object $\lnr_\A$. Moreover, since $M$-parametrized linear endomorphisms turn out to be simply linear endomorphisms in $\Simple_M[\X,\TT]$, each set $\lnr_\A(M)$ comes equipped with a rig structure that turns the functor $\lnr_\A(-)$ into a rig object in the presheaf category of $\X$. Since the presheaf category is always complete, one can construct the maximal subring of $\lnr_\A(-)$ in the presheaf category by taking the equalizer of the sum $+$ with the zero map $!0$ of $\lnr_\A(-)$. Thus, $\lnr_\A^-(-)$ is precisely the resulting ring object in the presheaf category.\\

The Lie functor $\gl(-)$ of the general linear group $\GL_\A$ of a differential object $\A$ is related to the functor $\lnr_\A(-)$ via a natural Lie algebra homomorphism, as proved by the following proposition.

\begin{proposition}
\label{proposition:the-map-chi}
Suppose that $\A$ is a differential object $\A$ which admits the general linear group, and that the tangent bundle $\T\GL_\A$ admits negatives\footnote{Negatives are required to define the Lie functor a group object. Notice that the tangent bundle of a Lie group has always negatives.}. Therefore, there is a Lie algebra homomorphism
\begin{align*}
&\chi_M^\A\colon\gl_\A(M)\to\lnr_\A^-(M)
\end{align*}
natural in $M$ which sends each $M$-\PLIVF $X\colon\GL_\A\times M\to\T\GL_\A$ to:
\begin{align*}
&\varphi_X\colon \A\times M\xrightarrow{\id\times\e\times\id}\A\times\GL_\A\times M\xrightarrow{z\times X}\T(\A\times\GL_\A)\xrightarrow{\T\Ev}\T\A\xrightarrow{\hat p}\A
\end{align*}
\end{proposition}
\begin{proof}
For starters, we show that $\varphi_X$ is an $M$-parametrized linear endomorphism of $\A$:
\begin{align*}
\varphi_X\lambda_\A&=~(\id\times\e\times\id)(z\times X)\T\Ev\hat p\lambda_\A    \tag{Definition of $\varphi_X$}\\
&=~(\id\times\e\times\id)(z\times X)\T\Ev\T\lambda_\A c\T\hat p         \tag{$\hat p\lambda_\A=\T\lambda_\A c\T\hat p$}\\
&=~(\id\times\e\times\id)(z\T\lambda_\A\times X\T z)\T\T\Ev c\T\hat p     \tag{Linearity of $\Ev$}\\
&=~(\id\times\e\times\id)(\lambda_\A z\times z\T X)\T\T\Ev\T\hat p      \tag{Naturality of $c$ and $z$}\\
&=~(\lambda_\A\times z)\T((\id\times\e\times\id)(z\times X)\T\Ev)\T\hat p\\
&=~(\lambda_\A\times z)\T\varphi_X
\end{align*}
Thus, we constructed a function $\chi_M^\A\colon\gl_\A(M)\to\lnr_\A(M)$ from the Lie functor to $\lnr_\A(M)$, which is the set of $M$-parametrized linear endomorphisms of $\A$. However, since, by assumption, the tangent bundle of $\GL_\A$ admits negatives, such a function factors through the natural inclusion $\lnr_\A^-(M)\to\lnr_\A(M)$. We shall now prove that this function preserves the Lie algebra structure. We start by showing that $\varphi_0$ is the zero morphism, where $0\=\pi_1z\colon\GL_\A\times M\to\T\GL_\A$:
\begin{align*}
\varphi_0&=~(\id\times\e\times\id)(z\times(\pi_1z)\T\Ev\hat p  \tag{Definition of $\varphi_0$}\\
&=~\pi_1(\id\times\e)z\T\Ev\hat p                  \tag{Naturality of $\pi_1$}\\
&=~\pi_1(\id\times\e)\Ev z\hat p                                 \tag{Naturality of $z$}\\
&=~\pi_1(\id\times\e)\Ev!0                                      \tag{$z\hat p=!0$}\\
&=~!_{\id\times M}0\\
&=~0
\end{align*}
Now, let us show that $\chi$ preserves the sum. Given two \PLIVFs $X,Y\colon\GL_\A\times M\to\T\GL_\A$:
\begin{align*}
\varphi_{X+Y}&=~(\id\times\e\times\id)(z\times(\<X,Y\>s))\T\Ev\hat p      \tag{Definition of $\varphi_{X+Y}$}\\      
&=~(\id\times\e\times\id)(\<z\times X,z\times Y\>)\T_2\Ev s\hat p  \tag{Linearity of $\Ev$ implies additivity of $\Ev$}\\
&=~(\id\times\e\times\id)(\<z\times X,z\times Y\>)\T_2\Ev\T_2\hat p+  \tag{$s\hat p=\T_2\hat p+$}\\
&=~\left\<(\id\times\e\times\id)(z\times X)\T\Ev\T\hat p,(\id\times\e\times\id)(z\times Y)\T\Ev\T\hat p\right\>+\\
&=~\varphi_X+\varphi_Y
\end{align*}
Finally, we shall prove that $\varphi_{[X,Y]}=[\varphi_X,\varphi_Y]$. Let us begin with considering the composition $\varphi_Y;_\KL\varphi_X$, where composition is in the coKleisli category $\Simple_M(\X,\TT)$:
\begin{align*}
\varphi_Y;_\KL\varphi_X&(\id\times\Delta_M)(\varphi_Y\times\id)\varphi_X
\end{align*}
After some basic rewriting and using the naturality of $z$, we can rewrite this as:
\begin{align*}
&\left((z\times\T\e\times z)(\T z\times\T Y)\T\T\Ev\T\hat p\times(\e\times\id)X\right)\T\Ev\hat p
\end{align*}
However, since $\Ev$ is linear in the $\A$-entry, $(\hat p\times\id)\Ev=(\id\times z)\T\Ev\hat p$. Using this identity, we can rewrite $\varphi_Y;_\KL\varphi_X$ as follows:
\begin{align*}
&\left((z\times\T\e\times z)(\T z\times\T Y)\T\T\Ev\times(\e\times\id)X\T z\right)\T\T\Ev\T\hat p\hat p
\end{align*}
However, using the left-invariance of $Y$, we can write the following commutative diagram:
\begin{equation*}
\begin{tikzcd}
M & \\
{\GL_\A\times M\times\GL_\A\times M\times M} \\
{\T(\GL_\A\times\GL_\A\times M)} & {\T(\GL_\A\times M)} \\
{\T\T(\GL_\A\times\GL_\A)} & {\T\T\GL_\A}
\arrow["{\<\e,\id,\e,\id,\id\>}"', from=1-1, to=2-1]
\arrow["{\<\T\e,z\>}", curve={height=-24pt}, from=1-1, to=3-2]
\arrow["{X\i\times X\times z}"', from=2-1, to=3-1]
\arrow["{\T(\m\times\id)}", from=3-1, to=3-2]
\arrow["{\T(z\times Y)}"', from=3-1, to=4-1]
\arrow["{\T Y}", from=3-2, to=4-2]
\arrow["{\T\T\m}"', from=4-1, to=4-2]
\end{tikzcd}
\end{equation*}
Furthermore, by construction, the multiplication $\m$ of $\GL_\A$ is uniquely defined by the following identity, $(\id\times\m)\Ev=(\Ev\times\id)\Ev$. Using this identity and the left-invariance of $Y$, we can rewrite $\varphi_Y;_\KL \varphi_X$ as follows:
\begin{align*}
&(z\times\e\times\id\times\e\times\id\times\e\times\id)(z\times((X\T z\T\i)\times((X\times z)\T Y)\times(X\T z)))\\
&\qquad(\T z\times\T\T\m_3)\T\T\Ev\T\hat p\hat p
\end{align*}
where $\m_3\=(\id\times\m)\m=(\m\times\id)\m$. In this expression there are three terms, $(X\T z\T\i)$, $((X\times z)\T Y)$, and $(X\T z)$, multiplied together by $\T\T\m_3$. However, the first and last terms invert each other. Furthermore, by an Eckmann-Hilton argument, the derivative of the multiplication $\m$ coincides with the sum of tangent vectors and thus is commutative. This allows us to swap the terms and simplify the former with the latter. What is left is the following expression:
\begin{align*}
&\varphi_Y;_\KL \varphi_X=((\id\times\e\times\id)(z\times X)\times z)(\T z\times\T Y)\T\T\Ev\T\hat p\hat p\\
&\qquad=
(\id\times\e\times\id)(z\T z\times(X;_M\T_MY))\T\T\Ev\T\hat p\hat p
\end{align*}
where $X;_M\T_MY$ denotes the composition of $X\colon\GL_\A\nto\T\GL_\A$ with $\T_MY\colon\T\GL_\A\nto\T\T\GL_\A$ in $\Simple_M(\X,\TT)$. Similarly, we also obtain that:
\begin{align*}
&\varphi_X;_\KL \varphi_Y=(\id\times\e\times\id)(z\T z\times(Y;_\KL\T^\KL X))\T\T\Ev\T\hat p\hat p
\end{align*}
Since $\A$ is a differential object, $\T\T\A\cong\A\times\A\times\A\times\A$. Thus, the canonical flip $c$ swaps the internal components, while $\T\hat p\hat p\colon\T\T\A\to\A$ projects the fourth component down to $\A$. Therefore, $\T\hat p\hat p=c\T\hat p\hat p$ (\cite[Proposition~3.6]{cockett:differential-bundles}). So, we can also rewrite:
\begin{align*}
&\varphi_X;_\KL \varphi_Y=(\id\times\e\times\id)(z\T z\times(Y;_\KL\T^\KL X)c)\T\T\Ev\T\hat p\hat p
\end{align*}
where we used the naturality of $c$ and that $z\T z c=z z=z\T z$. Now, let us unwrap the definition of $\varphi_{[X,Y]}$:
\begin{align*}
\varphi_{[X,Y]}&=~(\id\times\e\times\id)(z\times[X,Y])\T\Ev\hat p\\
&=~(\id\times\e\times\id)\left\{\left(\left(z\T z\times(X;_\KL\T^\KL Y)\right)-\left(z\T z\times(Y;_\KL\T^\KL Xc)\right)\right)\T\T\Ev\right\}\hat p
\end{align*}
By~\cite[Example~3.10]{cockett:connections}, if $\A$ is a differential object, $\<p_{\T \A}p_\A,T\hat p\hat p\>\colon\T\T\A\to\T\A$ defines a vertical affine connection of $\A$\footnote{Notice that in~\cite[Example~3.10]{cockett:connections}, the first component of $\T\A\cong\A\times\A$ contain the tangent vectors and not the points. We use the opposite convention.}. However, by~\cite[Lemma~3.13]{cockett:connections}, if $K\colon\T\T M\to\T M$ is an affine connecection, then $\{\varphi\}=\varphi K$. However, we can also compute:
\begin{align*}
&\{\varphi\}\hat p=\varphi\<p_{\T\A}p_\A,\T\hat p\hat p\>\hat p=\varphi\T\hat p\hat p
\end{align*}
Therefore, we can rewrite $\varphi_{[X,Y]}$ as:
\begin{align*}
&\left((\id\times\e\times\id)\left(z\T z\times(X;_\KL\T^\KL Y)\right)-(\id\times\e\times\id)\left(z\T z\times(Y;_\KL\T^\KL Xc)\right)\right)\\
&\qquad\T\T\Ev\T\hat p\hat p
\end{align*}
Finally, using the additivity of $\hat p$, we conclude that $\varphi_{[X,Y]}=\varphi_Y\varphi_X-\varphi_X\varphi_Y=[\varphi_X,\varphi_Y]$.
\end{proof}

The Lie algebra morphism $\chi\colon\gl_\A(M)\to\lnr_\A^-(M)$ compares the Lie functor of the general linear group $\GL_\A$ with the functor $\lnr_\A^-(-)$ that classifies linear endomorphisms of $\A$. When $\A$ admits the LEC object and its maximal subring $\lnr_\A^-$, then the functor $\lnr_\A^-(-)$ becomes representable. Furthermore, if also the general linear group $\GL_\A$ of $\A$ is a Lie group, then its Lie functor $\gl_\A(-)$ is also representable. We shall introduce a useful definition.

\begin{definition}
\label{definition:general-linear-lie-group}
We shall say that a differential object $\A$ admits the \textbf{general linear Lie group} provided that the general linear group object $\GL_\A$ of $\A$ exists and is a Lie group. We shall denote the differential Lie algebra of $\GL_\A$ by $\gl_\A$.
\end{definition}

We shall now prove that the comparison natural transformation $\chi$ is completely determined by the inclusion map $\theta_\A\colon\GL_\A\to\lnr_\A$, uniquely defined by $(\id\times\theta)\ev=\Ev$. To this end, we shall denote by $\rho\colon\lnr_\A^-\to\lnr_\A$ the inclusion map. Recall also that $\gl_\A\cong\T_\e\GL_\A$ and that $\lnr_\A\cong\T_\e\lnr_\A$, since $\lnr_\A$ is a differential object.

\begin{lemma}
\label{lemma:chi-vs-theta-maps}
If $\A$ admits the LEC object $\lnr_\A$, its maximal subring $\lnr_\A^-$, and the general linear Lie group $\GL_\A$, then there exists a linear Lie algebra homomorphism
\begin{align*}
&\chi_\A\colon\gl_\A\to\lnr_\A^-
\end{align*}
induced by the natural transformation $\chi_\A^M\colon\gl_\A(M)\to\lnr_\A(M)$. Furthermore:
\begin{align*}
&\chi_\A\rho_\A=\T_\e\theta_\A
\end{align*}
\end{lemma}
\begin{proof}
By the Yoneda lemma, there is a unique morphism $\chi'_\A\colon\gl_\A\to\lnr_\A$, which renders the following diagram
\begin{equation*}
\begin{tikzcd}
{\A\times\gl_\A} & \\
{\A\times\lnr_\A} & \A
\arrow["{\id\times\chi_\A'}"', from=1-1, to=2-1]
\arrow["{\varphi_\omega}", from=1-1, to=2-2]
\arrow["{\ev}"', from=2-1, to=2-2]
\end{tikzcd}
\end{equation*}
commutative, where $\varphi_\omega=\chi_{\gl_\A}^\A(\omega)$ and where $\omega\colon\GL_\A\times\gl_\A\to\T\GL_\A$ is the universal \PLIVF of $\GL_\A$. Let us unwrap the definition of $\varphi_\omega$:
\begin{align*}
\varphi_\omega&=~(\id\times\e\times\id)(z\times\omega)\T\Ev\hat p       \tag{Definition of $\varphi_\omega$}\\
&=~(\id\times\e\times\id)(z\times(z\times\iota_\e)\T\m)\T\Ev\hat p   \tag{Lemma~\ref{lemma:omega-vs-iota}}\\
&=~(\id\times\e\times\id)((z\times z)\T\Ev\times\iota_\e)\T\Ev\hat p   \tag{Definition of $\m$}\\
&=~((\id\times\e)\Ev z\times\iota_\e)\T\Ev\hat p      \tag{Naturality of $z$}\\
&=~(z\times\iota_\e)\T\Ev\hat p                          \tag{Definition of $\e$}
\end{align*}
Now, let us consider $\T_\e\theta_\A$, which is defined as the unique morphism, such that $\T_\e\theta_\A\<!\e,\lambda_{\lnr_\A}\>=\iota_\e\T\theta_\A$. Thus, by postcomposing on both sides by $\hat p_{\lnr_\A}$, we obtain that, $\T_\e\theta_\A=\T_\e\theta_\A\<!\e,\lambda_{\lnr_\A}\>\hat p_{\lnr_\A}=\iota_\e\T\theta_\A\hat p_{\lnr_\A}$, where we used that $\lambda_{\lnr_\A}\hat p_{\lnr_\A}$ is the identity. Thus:
\begin{align*}
(\id\times\T_\e\theta_\A)\ev&=~(\id\times(\iota_\e\T\theta_\A\hat p_{\lnr_\A}))\ev\\
&=~(\id\times\iota_\e\T\theta_\A)(\id\times\hat p_{\lnr_\A})\ev\\
&=~(\id\times\iota_\e\T\theta_\A)(z\times\id)\T\ev\hat p     \tag{Bilinearity}\\
&=~(z\times\iota_\e)\T((\theta_\A\times\id)\ev)\hat p\\
&=~(z\times\iota_\e)\T\Ev\hat p                                          \tag{Definition of $\theta_\A$: $(\theta_\A\times\id)\ev=\Ev$}
\end{align*}
Using the universal property of $\ev$, we conclude that $\chi_\A'=\T_\e\theta_\A$. However, since every morphism $\T_\e f$ is always linear, $\chi_\A'$ is also linear. To conclude, it is only left to show that $\chi_\A'=\chi_\A\rho$ for some $\chi_\A\colon\gl_\A\to\lnr_\A^-$. However, since $\gl_\A$ has negatives given by $\T_\e\i$, where $\i$ denotes the inverse map of $\GL_\A$, it follows that $\chi_\A'$ lands inside of $\lnr_\A^-$, that is, it factors through $\rho\colon\lnr_\A^-\to\lnr_\A$. Using the universal property of $\rho$ that establishes that $\<\rho,-\rho\>$ is an equalizer map of $+$ and $!0$ as in Definition~\ref{definition:maximal-differential-subring}, one can easily prove that $\chi_\A$ is in fact a linear Lie algebra homomorphism.
\end{proof}

Before proving the main result, we indulge in another technical lemma.

\begin{lemma}
\label{lemma:multiplication-and-unit-at-e}
By identifying the tangent space $\T_\e\lnr_\A$ of $\lnr_\A$ at $\e\colon\*\to\lnr_\A$ with $\lnr_\A$ itself, the morphism $\T_\e\m$ coincides with the sum $+$ of $\lnr_\A$ and $\T_\e\e$ with the zero $0$.
\end{lemma}
\begin{proof}
To prove this lemma, one can argue as in~\cite[Theorem~2]{cockett:lie-groups-tangent-cats}, in which Cockett and Schwarz used an Eckmann-Hilton argument to prove that $\T_\e\m=+$ and $\T_\e\e=0$, for Lie groups.
\end{proof}

Thanks to Theorem~\ref{theorem:triple-equalizer}, the general linear group $\GL_\A$ is the tangent triple equalizer of $\m$, $!\e$, and $\m^\op$ of $\lnr_\A$. Using that result, we now prove the main theorem of this paper: that the differential Lie algebra of $\GL_\A$ is isomorphic to the differential Lie algebra $\lnr_\A^-$.

\begin{theorem}
\label{theorem:main-theorem}
Suppose that $\A$ admits the general linear Lie group object $\GL_\A$ and the LEC object $\lnr_\A$. Then, the maximal subring $\lnr_\A^-$ of $\lnr_\A$ exists, and the differential Lie algebra $\gl_\A$ of $\GL_\A$ is isomorphic to the differential Lie algebra $\lnr_\A^-$. In particular, the map $\chi$ of Lemma~\ref{lemma:chi-vs-theta-maps} is an isomorphism of differential Lie algebras.
\end{theorem}
\begin{proof}
From Theorem~\ref{theorem:triple-equalizer}, $\GL_\A$ can be realized as the tangent triple equalizer of $\m$, $!\e$, and $\m^\op$ of $\lnr_\A$. By applying $\T_\e$, which takes the tangent space at $\e$ to this tangent triple equalizer, we obtain a new tangent triple equalizer:
\begin{equation*}
\begin{tikzcd}
{\T_\e\GL_\A} & {\T_\e\lnr_\A\times\T_\e\lnr_\A} & {\*} & {\T_\e\lnr_\A}
\arrow[from=1-1, to=1-2]
\arrow[from=1-2, to=1-3]
\arrow["{\T_\e\m}", curve={height=-24pt}, from=1-2, to=1-4]
\arrow["{\T_\e\m^\op}"', curve={height=24pt}, from=1-2, to=1-4]
\arrow["{\T_\e\e}", from=1-3, to=1-4]
\end{tikzcd}
\end{equation*}
In fact, $\T_\e$ is constructed by means of tangent pullbacks; however, as limits commute with other limits, tangent limits commute with other tangent limits. Furthermore, since $\lnr_\A$ is a differential object, $\T_\e\lnr_\A\cong\lnr_\A$. However, by Lemma~\ref{lemma:multiplication-and-unit-at-e}, $\T_\e\m=+$ and $\T_\e\e=0$. Moreover, since $+$ is commutative, it also follows that $\T_\e\m^\op=+^\op=+$. Therefore, we can rewrite the tangent triple equalizer diagram as the following tangent equalizer diagram:
\begin{equation*}
\begin{tikzcd}
{\gl_\A} & {\lnr_\A\times\lnr_\A} && {\lnr_\A}
\arrow[from=1-1, to=1-2]
\arrow["{+}", shift left=2, from=1-2, to=1-4]
\arrow["{!0}"', shift right=2, from=1-2, to=1-4]
\end{tikzcd}
\end{equation*}
However, by Definition~\ref{definition:maximal-differential-subring}, the tangent equalizer of $+$ and $!0$ is precisely the maximal differential subring $\lnr_\A^-$ of $\lnr_\A$. To conclude, we need to prove that the unique isomorphism $\gl_\A\cong\lnr_\A^-$ must be $\chi_\A$. This will prove that not only $\gl_\A$ and $\lnr_\A^-$ are isomorphic as differential objects but that they are also isomorphic as differential Lie algebras, since by Lemma~\ref{lemma:chi-vs-theta-maps}, $\chi$ is a Lie algebra homomorphism. To prove this statement, recall that $\chi\rho=\T_\e\theta$ as established by Lemma~\ref{lemma:chi-vs-theta-maps}, where $\<\rho,-\rho\>\colon\lnr_\A^-\to\lnr_\A\times\lnr_\A$ is the equalizer map of $+$ and $!0$. Let $\xi\colon\gl_\A\to\lnr_\A^-$ denote the unique isomorphism that satisfies $\xi\<\rho,-\rho\>=\T_\e\varphi=\<\T_\e\theta,\T_\e\i\T_\e\theta\>$. By postcomposing on both sides by $\pi_1$, we obtain that $\xi\rho=\T_\e\theta=\chi\rho$. However, since $\rho$ is monic, this implies that $\xi=\chi$. Therefore, $\chi$ becomes a differential Lie algebra isomorphism.
\end{proof}

When the differential object $\A$ admits negatives, then $\lnr_\A$ also admits negatives, with $-$ defined as the unique morphism that makes the following diagram
\begin{equation*}
\begin{tikzcd}
{\A\times\lnr_\A} & \A \\
{\A\times\lnr_\A} & \A
\arrow["\ev", from=1-1, to=1-2]
\arrow["{\id\times-}", from=2-1, to=1-1]
\arrow["\ev"', from=2-1, to=2-2]
\arrow["{-}"', from=2-2, to=1-2]
\end{tikzcd}
\end{equation*}
commutative. Therefore, when $\A$ admits negatives, $\lnr_\A\cong\lnr_\A^-$. In this case, Theorem~\ref{theorem:main-theorem} becomes the classic result of differential and algebraic geometry that establishes that the Lie algebra of the general linear group $\GL_n(\Real)$ is isomorphic to the Lie algebra $\lnr_n(\Real)$ of linear endomorphisms of $\Real^n$.

In particular, thanks to~\cite[Lemma~2.2.1 and Corollary~2.2.5]{macadam:vector-bundles}, when the tangent structure comes with negatives, in that the tangent bundle $\T M\to M$ of each object admits a negation $n\colon\T M\to\T M$, every differential bundle also comes with negatives. So, in particular, each differential object admits negatives.

\begin{corollary}
\label{corollary:main}
If $\A$ admits negatives, the Lie algebra of the general linear Lie group $\GL_\A$ of $\A$ is isomorphic to the Lie algebra $\lnr_\A\cong\lnr_\A^-$. In particular, in a Cartesian tangent category with negatives, $\gl_\A\cong\lnr_\A$.
\end{corollary}


\printbibliography

\end{document}